\documentclass[reqno]{amsart}
\usepackage{amsmath, amsthm, amssymb, xypic}
\usepackage{hyperref}
\usepackage[nobysame, non-sorted-cites, alphabetic]{amsrefs}
\usepackage{enumitem}
\setenumerate{label=\textnormal{(\arabic*)}} 
\usepackage{graphicx}

\usepackage{color}

\xyoption{rotate}

\newcommand{\cat}[1]{\ensuremath{\mathsf{#1}}}
\newcommand{\gen}[1]{\ensuremath{\langle #1 \rangle}}

\newcommand{\op}{\ensuremath{^\mathrm{op}}}

\newcommand{\ctens}{\widehat{\otimes}}

\DeclareMathOperator{\Ring}{\cat{Ring}}
\newcommand{\Ringop}{\Ring\op}
\DeclareMathOperator{\cRing}{\cat{cRing}}
\newcommand{\cRingop}{\cRing\op}
\DeclareMathOperator{\Set}{\cat{Set}}

\DeclareMathOperator{\Top}{\cat{Top}}

\DeclareMathOperator{\Alg}{\cat{Alg}}
\DeclareMathOperator{\cAlg}{\cat{cAlg}}
\DeclareMathOperator{\fdAlg}{\cat{fdAlg}}
\DeclareMathOperator{\fdcAlg}{\cat{fdcAlg}}
\DeclareMathOperator{\PCAlg}{\cat{PCAlg}}

\DeclareMathOperator{\Coalg}{\cat{Coalg}}
\DeclareMathOperator{\fdCoalg}{\cat{fdCoalg}}
\DeclareMathOperator{\fdcCoalg}{\cat{fdcCoalg}}

\DeclareMathOperator{\Sch}{\cat{Sch}}
\DeclareMathOperator{\fSch}{\cat{fSch}}

\DeclareMathOperator{\fin}{fin}

\DeclareMathOperator{\Spec}{Spec}
\renewcommand{\O}{\mathcal{O}}
\DeclareMathOperator{\Max}{Max}

\DeclareMathOperator{\Hom}{Hom}
\DeclareMathOperator{\End}{End}

\DeclareMathOperator{\ann}{ann}
\newcommand{\m}{\mathfrak{m}}

\DeclareMathOperator{\pts}{pts}
\DeclareMathOperator{\Dist}{Dist} 
\DeclareMathOperator{\corad}{corad}

\DeclareMathOperator{\Br}{Br}
\DeclareMathOperator{\Azu}{\mathcal{A}}

\newcommand{\F}{\mathcal{F}}
\newcommand{\eps}{\varepsilon}
\newcommand{\dual}{\widehat}

\newcommand{\invlim}{\varprojlim}
\newcommand{\dirlim}{\varinjlim}

\newcommand{\Z}{\mathbb{Z}}

\newcommand{\C}{\mathbb{C}}
\newcommand{\M}{\mathbb{M}}

\newcommand{\A}{\mathbb{A}}
\newcommand{\T}{\mathbb{T}}

\newcommand{\catC}{\mathcal{C}}

\DeclareMathOperator{\lMod}{\!-\cat{Mod}}

\DeclareMathOperator{\lPC}{\!-\cat{PC}}

\newcommand{\e}{\'{e}}

\newcommand{\separate}{\bigskip}

\numberwithin{equation}{section}

\theoremstyle{plain}
\newtheorem{theorem}[equation]{Theorem}
\newtheorem{corollary}[equation]{Corollary}
\newtheorem{lemma}[equation]{Lemma}
\newtheorem{proposition}[equation]{Proposition}

\theoremstyle{definition}
\newtheorem{definition}[equation]{Definition}

\newtheorem{example}[equation]{Example}
 
\newtheorem{remark}[equation]{Remark}

\begin{document}
\title[The finite dual as a quantized maximal spectrum]{The finite dual coalgebra as a quantization of the maximal spectrum}
\author{Manuel L. Reyes}
\address{Department of Mathematics\\ University of California, Irvine\\
  419 Rowland Hall\\ Irvine, CA 92697--3875\\ USA}
\email{mreyes57@uci.edu}
\urladdr{https://www.math.uci.edu/~mreyes/}
\date{January 26, 2024}
\thanks{This work was supported in part by NSF grant DMS-2201273.}

\begin{abstract}
In pursuit of a noncommutative spectrum functor, we argue that the Heyneman-Sweedler finite dual coalgebra can be viewed as a quantization of the maximal spectrum of a commutative affine algebra,
integrating prior perspectives of Takeuchi, Batchelor, Kontsevich-Soibelman, and Le~Bruyn. 
We introduce \emph{fully residually finite-dimensional} algebras $A$ as those with enough finite-dimensional representations to let $A^\circ$ act as an appropriate depiction of the noncommutative maximal spectrum of $A$; importantly, this class includes affine noetherian PI algebras. 
In the case of prime affine algebras that are module-finite over their center, we describe how the Azumaya locus is represented in the finite dual.
This is used to describe the finite dual of quantum planes at roots of unity as an endeavor to visualize the noncommutative space on which these algebras act as functions. Finally, we discuss how a similar analysis can be carried out for other maximal orders over surfaces.
\end{abstract}

\maketitle

\setcounter{tocdepth}{2}
\makeatletter
\def\l@subsection{\@tocline{2}{0pt}{2.5pc}{5pc}{}}
\def\l@subsubsection{\@tocline{2}{0pt}{5pc}{7.5pc}{}} 
\makeatother
\tableofcontents

\section{Introduction}
\label{sec:intro}

\subsection{Prologue: noncommutative spectral theory}\label{sub:spectral}

This paper is intended as a contribution to \emph{noncommutative spectral theory}, by which we mean the development and study of invariants of associative rings and algebras that extend the usual spectrum of a commutative ring or algebra, with the goal of gaining insight into the structure of noncommutative rings. Ideally one would wish for a noncommutative spectrum that can be equipped with enough extra structure to form a complete invariant, allowing us to recover the ring up to isomorphism. The basis for this hope is the Zariski spectrum of a commutative ring $R$, which can be equipped with its structure sheaf $\O_{\Spec(R)}$ that allows us to recover $R$ up to isomorphism as the ring of global sections of the affine scheme $(\Spec(R), \O_{\Spec(R)})$.
Probably the best-known work on noncommutative spectral theory regards the development of structure sheaves on noncommutative prime spectra. Considerable efforts have been made in this direction over several decades; see~\cite{vanoystaeyen:sheaves} for a comprehensive survey and a large bibliography on this topic.

When considering extensions of the Zariski spectrum as an assignment from rings to topological spaces, we believe that there are two important properties that one should ask of such a spectrum: 
\begin{enumerate}
\item the assignment returns the usual Zariski spectrum of a commutative ring, 
\item the assignment extends to a functor $\Ringop \to \Top$. 
\end{enumerate}
Condition~(2) would bring us closer to a duality between noncommutative rings and appropriate spatial objects.
Aside from aesthetic considerations, it is also important because it guarantees that partial information about the spectrum of a ring can be obtained from the spectra of its subrings and quotient rings.
It is also crucial if one wishes for such a spectrum to take into account algebraic quantum theory~\cite{landsman:aqm}, where commutative subalgebras of a noncommutative algebra are intimately tied to the relationship between classical and quantum information~\cite{heunen:classicalfaces}. 

Unfortunately, the results of~\cites{reyes:obstructing, bmr:kochenspecker} indicate a major obstruction to this hope: any functor satisfying the above properties must assign the empty space to any ring with $n \times n$ matrix structure for $n \geq 3$. One might conclude from this obstruction that the problem somehow lies in the points themselves, particularly if one is familiar with point-free approaches to topology~\cite{johnstone:pointless} such as locales, toposes, and quantales. However, subsequent results~\cites{bergheunen:obstructions, reyes:sheaves, bmr:kochenspecker} indicate that similar obstructions persist at this level of generality.

In our view, these obstructions indicate that the usual building blocks of topology, whether from point sets or open covers without points, are too commutative in nature to serve as a deep spectral invariant for noncommutative rings. This highlights a fundamental question that has been minimized for too long: 
\begin{quote}
What objects should play the role of \emph{noncommutative sets}  within noncommutative geometry? 
\end{quote} 
While many kinds of noncommutative spaces have been studied across the various branches of noncommutative geometry, the most basic case of noncommutative \emph{discrete} spaces (i.e., sets) has received comparatively little attention. The obstructions discussed above suggest that we cannot make serious progress in noncommutative spectral theory until this gap in noncommutative mathematics has been filled.  Fortunately, recent work such as~\cite{heunenreyes:discretization, mrv:quantum functions, kornell:quantumsets, kornell:discrete, rump:quantum} has begun to address this issue, although it is limited to the setting of noncommutative geometry based on C*-algebras. 

In light of the above, we will use the term \emph{functorial spectral theory} to describe the pursuit of a functorial invariant of associative rings (or C*-algebras) that extends the Zariski (or Gelfand) spectrum from commutative to noncommutative rings. The major problem becomes locating a suitable category $\mathfrak{S}$ of \emph{noncommutative sets}, with a full and faithful functor $\Set \hookrightarrow \mathfrak{S}$, such that there is a ``quantum spectrum'' functor $\Sigma \colon \Ringop \to \mathfrak{S}$ making the following diagram commute (up to natural isomorphism): 
\begin{equation}\label{eq:Spec diagram}
\begin{split}
\xymatrix{
\cRingop \ar[r]^-{\Spec} \ar@{^{(}->}[d]  & \Top \ar[r]^U & \Set \ar@{^{(}->}[d] \\
\Ringop \ar[rr]^-{\Sigma} &  & \mathfrak{S}
}
\end{split}
\end{equation}
Ideally, one should be able to equip noncommutative sets with (some noncommutative version of) a topology in order to endow each quantum spectrum with a structure sheaf, from which one can recover the original ring up to isomorphism. While this is a daunting problem, there have been recent advances~\cites{heunenreyes:active, heunenreyes:discretization} in functorial spectral theory for C*-algebras that give us reason to maintain hope in the face of this challenge.

\subsection{Quantizing the maximal spectrum}\label{sub:quantum max}

A comprehensive solution to the problems raised above remains out of reach for now. In this paper we pursue the more modest goal of extending the \emph{maximal spectrum} functor from the category $\cat{cAff}$ of commutative affine algebras over a field~$k$ to a suitable class of algebras that are only ``mildly noncommutative.'' Our intent is to catch a glimpse of the underlying discrete objects of noncommutative spaces, even if we cannot directly access the objects themselves. 

Our point of departure lies in the algebraic approach to quantum mechanics~\cite{landsman:aqm}, where commutative and noncommutative algebras are respectively embodied as classical and quantum observables. We believe that there are deep reasons to attempt to build a functorial spectral theory that draws inspiration from quantum physics, in part because many of the obstructions mentioned above hinge upon the Kochen-Specker Theorem of quantum mechanics~\cite{bmr:kochenspecker} and in part due to the important role that commutative subalgebras of noncommutative algebras play in each theory~\cite{heunen:classicalfaces, reyes:obstructing}. 

Let us imagine a commutative algebra as an algebra of observables in classical mechanics, so that its spectrum corresponds to the space of states of a classical system. In passing to quantum mechanics, the \emph{superposition principle} forces physicists to allow for linear combinations (more precisely, convex combinations) of quantum states. Thus for a commutative $k$-algebra $A$ we aim to replace the set $\Max(A)$ with $k$-linear combinations of its points. (We may imagine that $k$ is algebraically closed for the moment, to minimize technicalities.) It turns out that we may recover any set from its $k$-linear span when endowed with just a bit of structure: a \emph{coalgebra}. This yields a fully faithful embedding of the category $\Set$ into the category $\Coalg$ of $k$-coalgebras (Subsection~\ref{sub:quantized sets}).
In this way we think of coalgebras as ``linear spans of quantum states,'' with comultiplication functioning as a ``quantum diagonal'' map. We can imagine elements of the coalgebra as linear combinations of points, but it is not possible in general to isolate individual points. Note that Batchelor similarly viewed coalgebras as generalized sets in noncommutative geometry in~\cite{batchelor:ncg}. 
At the same time, we describe in Subsection~\ref{sub:distributions} how each $k$-scheme has a cocommutative \emph{coalgebra of distributions} that can be viewed as underlying discrete objects. This idea dates back to work of Takeuchi~\cite{takeuchi:tangent, takeuchi:formal}, and we connect our point of view to the definitions used in that work. A similar perspective in differential geometry arose in work of Batchelor~\cite{batchelor:maps}.

The \emph{finite dual coalgebra} construction of Heyneman and Sweedler~\cite[1.3]{heynemansweedler} associates to every $k$-algebra $A$ a coalgebra $A^\circ$ (whose definition we recall in Subsection~\ref{sub:Sweedler}) in such a way that it forms a functor
\[
(-)^\circ \colon \Alg\op \to \Coalg.
\]
For commutative affine algebras $A$, the coalgebra $A^\circ$ agrees with the coalgebra of distributions on the scheme $\Spec(A)$, and its simple subcoalgebras are in bijection with $\Max(A)$. \
Thus for $k$-algebras that are endowed with ``many'' finite-dimensional representations, we argue that the
finite dual coalgebra serves as a suitable approximation to the quantized underlying discrete object of its noncommutative affine scheme.
More precisely, it can be viewed as containing a linearization of the maximal spectrum of an
affine scheme of finite type over $k$.
We point out that dual coalgebras were similarly viewed as distributions on noncommutative spaces in the noncommutative thin schemes~\cite[Section~2]{kontsevichsoibelman} of Kontsevich and Soibelman, and the perspective of dual coalgebras as noncommutative spectral objects has been explored further by Le Bruyn in~\cite{lebruyn:dual} and~\cite[Section~3.2]{lebruyn:branes}.

Embedding the category of sets into the category of coalgebras as described above, we obtain a diagram 
\begin{equation}\label{eq:Max diagram}
\begin{split}
\xymatrix{
\cat{cAff}\op \ar[r]^-{\Max} \ar@{^{(}->}[d]  & \Top \ar[r]^U & \Set \ar@{^{(}->}[d] \\
\Alg\op \ar[rr]^{(-)^\circ} &  & \mathbf{\Coalg}
}
\end{split}
\end{equation}
that is similar to~\eqref{eq:Spec diagram}, but which \emph{does not commute} up to isomorphism. Thus the dual coalgebra only approximates the space of closed points in a noncommutative spectrum. This is a necessary evil if we wish to require functoriality and nontrivial behavior for matrix algebras, as explained precisely in Theorem~\ref{thm:subfunctor}. 
While we remain cautiously optimistic that diagram~\eqref{eq:Spec diagram} might one day be realized, the situation in~\eqref{eq:Max diagram} can serve as a model for functorial invariants that do not strictly extend the Zariski spectrum, but rather extend some mild enhancement of it.

As alluded to above, it is only reasonable to expect that the dual coalgebra $A^\circ$ is a reasonable substitute for a maximal spectrum if $A$ has sufficiently many finite-dimensional representations. In Subsection~\ref{sub:fully RFD} we introduce a class of algebras $A$ for which we believe that $A^\circ$ captures a relatively full picture of the noncommutative geometry of $A$. These are called \emph{fully residually finite-dimensional (RFD)} algebras, and their defining property roughly states that the finite-dimensional representation theory of $A$ is rich enough to capture information about all finitely generated modules. Importantly, the fully RFD algebras contain a large class of familiar noncommutative algebras: the affine noetherian PI algebras. 

It is certainly true that the fully RFD property is quite restrictive as a ring-theoretic condition. However, a common phenomenon within noncommutative algebraic geometry is that many $q$-deformed families of quantum algebras contain affine  algebras that are module-finite over their center, typically when the deforming parameters are equal to roots of unity~\cite[Part~III]{browngoodearl:quantum}. Such an algebra satisfies a polynomial identity and is noetherian (see~\cite[Corollary~13.1.13, Lemma~13.9.10]{mcconnellrobson}). Thus these quantum families of algebras have nontrivial intersection with the class of fully RFD algebras, and one can try to obtain a geometric understanding of algebras in this intersection by studying the finite dual as a quantized maximal spectrum.

%

For affine algebras $A$ that are module-finite over their center $Z(A)$, functoriality guarantees that the dual coalgebras are related by a morphism $A^\circ \to Z(A)^\circ$. In this case, we aim to give a clearer relationship between $A^\circ$ and the geometry of $\Spec Z(A)$ in Section~\ref{sec:Azumaya}. 
Assuming that $A$ is prime, we describe $A^\circ$ in terms of a large part that corresponds to the Azumaya locus~\cite[Section~3]{browngoodearl:pi} tensored with a matrix coalgebra, with a direct sum complement that corresponds to the a noncommutative formal neighborhood of a closed subscheme of $\Spec Z(A)$. 
We use this to describe the dual coalgebra of the quantum plane $\O_q(k^2)$ where $k$ is algebraically closed and $q \in k^*$ is a root of unity, hoping to provide a glimpse of a concrete mathematical object that can be visualized as a quantum plane. 
This example underscores the potential of the perspective described above to place geometric and quantum language on equal footing. 
In the general case, we also discuss how knowledge of the local structure of $A$ relative to the maximal spectrum of $Z(A)$ can be used to describe $A^\circ$.
We hope that this will inspire future efforts to visualize noncommutative spaces in a similar manner.

\subsection{Outline of the paper}


Section~\ref{sec:discrete} begins by recalling the construction of the finite dual coalgebra. We then describe how coalgebras can be viewed both as quantized sets and as underlying discrete objects for (commutative) schemes of finite type over a field. 

This motivates the main thesis, described in Section~\ref{sec:quantized Max}, that the finite dual functor is a reasonable replacement for the maximal spectrum for algebras that are ``not too far'' from being commutative. Because $A^\circ$ is strictly larger than the linearized maximal spectrum for commutative affine algebras $A$, we take time to explain in Theorem~\ref{thm:subfunctor} that it is, in a sense, the best choice for a coalgebra-valued spectrum functor whose behavior on matrix algebras agrees with quantum intuition. We then introduce \emph{fully residually finite-dimensional} algebras as those algebras whose finite dual offers a useful substitute for the maximal spectrum.
Proposition~\ref{prop:FBN RFD} shows that this class of algebras includes all affine noetherian PI algebras, an important class of algebras that has a nontrivial intersection with many interesting and important quantized families of noncommutative algebras.


Finally, in Section~\ref{sec:Azumaya} we restrict to the case of affine algebras over an algebraically closed field that are module-finite over their center. For such algebras that are prime, we provide a direct sum decomposition of the dual coalgebra in which one of the summands corresponds to the Azumaya locus. This is applied to describe the dual coalgebra of the quantum plane $\O_q(k^2)$ for $k$ an algebraically closed field and $q \in k$ a root of unity. The paper concludes with a discussion of how similar analyses can be carried out for orders over affine $k$-varieties.

\section{Coalgebras as noncommutative discrete spaces}\label{sec:discrete}

In this section we will describe how coalgebras can be viewed as discrete geometric objects. 
This begins with a brief recollection of dual coalgebras in Subsection~\ref{sub:Sweedler}.
We discuss coalgebras as quantized sets in Subsection~\ref{sub:quantized sets}, where we attempt to provide some more intuition behind this viewpoint. Then in Subsection~\ref{sub:distributions} we explain how cocommutative coalgebras arise in algebraic geometry as collections of distributions on $k$-schemes.
Taken together, these suggest how coalgebras serve as discrete objects in noncommutative geometry.

Let $k$ be a field, which is completely arbitrary except where explicitly stated otherwise. Unadorned tensor symbols $- \otimes -$ denote tensor over $k$. In this paper, all algebras are unital and associative, all coalgebras are counital and coassociative, and morphisms of these objects preserve the (co)unit. 
We let $\Alg = \Alg_k$ denote the category of $k$-algebras and $\Coalg = \Coalg_k$ denote the category of $k$-coalgebras.

\subsection{Reminder on dual coalgebras}\label{sub:Sweedler}

Let $A$ be a finite-dimensional $k$-algebra. 
The natural map
\[
A^* \otimes A^* \to (A \otimes A)^*,
\]  
defined by allowing pure tensors $\phi \otimes \psi \in A^* \otimes A^*$ to act as functionals on $a \otimes b \in A \otimes A$ via $(\phi \otimes \psi)(a \otimes b) = \phi(a)\psi(b)$, is an isomorphism because $A$ is finite-dimensional. Thus the multiplication $m \colon A \otimes A \to A$ and unit $\eta \colon k \to A$ maps respectively dualize to a comultiplication and counit
\begin{align*}
\Delta &= m^* \colon A^* \to (A \otimes A)^* \cong A^* \otimes A^*, \\
\epsilon &= \eta^* \colon A^* \to k^* = k.
\end{align*}
These satisfy coassociativity and counitality as duals of associativity and unitality, so that $(A^*, \Delta, \epsilon)$ is a coalgebra.

%
%
\separate

For an infinite-dimensional algebra $A$, the natural embedding $A^* \otimes A^* \to (A \otimes A)^*$ is not an isomorphism, so the full linear dual $A^*$ does not naturally inherit a coalgebra structure. Nevertheless, there is a subspace $A^\circ \subseteq A^*$ that \emph{does} naturally form a coalgebra, called the \emph{finite dual} (occasionally called the \emph{Sweedler dual}).
For details regarding the following discussion, see~\cite[Chapter~VI]{sweedler:hopf} or~\cite[Section~1.5]{dnr:hopf}.

Letting $m \colon A \otimes A \to A$ denote the multiplication of $A$,
the subspace $A^\circ \subseteq A^*$ is defined to be the set of those $\phi \in A^*$ that satisfy
the following equivalent conditions:
\begin{itemize}
\item[(SD1)] $m^*(\phi) \in (A \otimes A)^*$ lies in the subspace $A^* \otimes A^* \subseteq (A \otimes A)^*$;
\item[(SD2)] the kernel of $\phi$ contains an ideal $I \unlhd A$ of finite codimension in $A$.
\end{itemize}
It turns out that the restriction of $m^*$ to $A^\circ$ has image in $A^\circ \otimes A^\circ
\subseteq A^* \otimes A^*$. In this way, $\Delta = m^*$ restricts to a comultiplication 
\[
\Delta \colon A^\circ \to A^\circ \otimes A^\circ
\] 
making the finite dual a coalgebra, whose counit is the dual of the unit of $A$.

The coalgebra structure on $A^\circ$ can be alternately described in terms of condition~(SD2) above.
Within the lattice of ideals of $A$, the family 
\begin{equation}\label{eq:F(A)}
\F(A) = \{I \unlhd A \mid \dim_k(A/I) < \infty\}
\end{equation}
of ideals having finite codimension forms a filter: it is upward-closed and closed under pairwise intersections. Thus the 
diagram of finite-dimensional algebras of the form $A/I$ for $I \in \F(A)$ forms an inversely directed
system. As each $\phi \in A^\circ$ is induced by some $\overline{\phi} \in (A/I)^*$ via the canonical
homomorphism $A \to A/I$ for some $I \in \F(A)$, we have an isomorphism of vector spaces
\begin{equation}\label{eq:Sweedler colimit}
A^\circ \cong \dirlim_{I \in \F(A)} (A/I)^*.
\end{equation}
Since the algebras $A/I$ above are finite-dimensional, the above is a directed limit of finite-dimensional
coalgebras, and the isomorphism~\eqref{eq:Sweedler colimit} is in fact an isomorphism of coalgebras.

Let $\Alg$ denote the category of $k$-algebras and their homomorphisms.
One can check using either condition~(SD1) or~(SD2) that any algebra homomorphism $\phi \colon A \to B$ dualizes to a morphism of coalgebras $B^\circ \to A^\circ$, so that the finite dual forms a functor
\begin{equation}\label{eq:Sweedler}
(-)^\circ \colon \Alg\op \to \Coalg.
\end{equation}
In fact, it enjoys the following adjoint relationship with the dual algebra functor $(-)^* \colon \Coalg\op \to \Alg$: for any algebra $A$ and coalgebra $C$, one has natural isomorphisms
\begin{equation}\label{eq:adjoint}
\Alg(A,C^*) \cong \Coalg(C, A^\circ).
\end{equation}

%

%

\subsection{Coalgebras as generalized sets}\label{sub:quantized sets}

As stated in the introduction, we will view coalgebras as the underlying discrete objects of noncommutative
spaces. To motivate this perspective, we will explain in detail how sets can be ``linearized'' in order
to provide a full and faithful embedding into the category of $k$-coalgebras. 


For a set $X$, we let $kX$ denote the $k$-vector space with basis $X$, which we call the
\emph{linearization} of $X$. 
%
This forms a $k$-coalgebra $(kX, \Delta, \epsilon)$ if we endow it with the comulitplication and counit
\begin{align*}
\Delta &\colon kX \to kX \otimes kX, \\
\epsilon &\colon kX \to k,
\end{align*}
that are defined on each basis element $x \in X$ by $\Delta(x) = x \otimes x$ and $\epsilon(x) = 1$. 
(This construction is well known, as in~\cite[Example~1.1.4]{dnr:hopf} 
and~\cite[III.1, Example~3]{kassel:quantum}.)
It is evident that linearization forms a functor
\[
k- \colon \Set \to \Coalg
\]
that is (unsurprisingly) defined on morphisms by linear extension.

This linearization functor has a right adjoint, constructed as follows. 
For a $k$-coalgebra $(Q,\Delta,\epsilon)$, we say that a nonzero element $x \in Q$ is \emph{pointlike}\footnote{In the literature, such an element is typically called \emph{grouplike}. We use this alternate terminology in order to emphasize a geometric perspective on coalgebras, especially ones that do not arise from any connection to a Hopf algebra or bialgebra. Fortunately, this slightly modified term fits well with the established notion of \emph{pointed} coalgebras as well as \emph{points} of an object in a category as global elements.} 
(or is a \emph{$k$-point}) if $\Delta(x) = x \otimes x$; the counit axiom yields $x = \epsilon(x)x$, so 
that $x \neq 0$ ensures $\epsilon(x) = 1$.
We let $\pts(Q)$ denote the set of $k$-points of any coalgebra $Q$. As any coalgebra morphism 
preserves pointlike elements, the assignment of $k$-point-sets forms a functor 
$\pts \colon \Coalg \to \Set$, which is easily seen to be naturally isomorphic to the functor 
$\Coalg(k,-)$. 
Now for a set $X$ and a coalgebra $Q$, a routine verification yields the following adjoint isomorphism:
\[
\Coalg(kX,Q) \cong \Set(X,\pts(Q)).
\]
As the terminology suggests, a coalgebra $kX$ defined from a set $X$ satisfies $\pts(kX) = X$. 
Then for sets $X$ and $Y$, the adjunction above restricts to
\[
\Coalg(kX,kY) \cong \Set(X,\pts(kY)) = \Set(X,Y),
\]
showing that linearization is fully faithful.

In this way we can view sets as a full subcategory of $\Coalg$ after linearization, allowing us to treat coalgebras as generalizations of (linearized) sets.
As mentioned in Section~\ref{sec:intro}, we in fact see coalgebras as a \emph{quantization} of sets, where $k$-linear combinations of points are imagined as a \emph{superposition} (over $k$) of states. The preceding discussion illustrates that the added structure of a coalgebra allows us to recover the individual points of $X$ from $kX$. 

\begin{example}\label{ex:coalgebra examples}
We briefly recall a few examples of well-known coalgebras, to which we will refer later.
\begin{enumerate}
\item \label{ex:matrix coalgebra} The \emph{comatrix coalgebra} $\M^d = \M^d(k)$ has basis $\{E^{ij} \mid 1 \leq i,j \leq d\}$ with
\[
\Delta(E^{ij}) = \sum\nolimits_r E^{ir} \otimes E^{rj} \quad \mbox{and} \quad
\epsilon(E^{ij}) = \delta_{ij}.
\]
Identifying the $E^{ij}$ with the dual basis of the matrix units $E_{ij} \in \M_d(k)$, we have $\M^d \cong (\M_d(k))^*$. 
By this duality, we view $\M^d$ as a $d$-level (or $d$-state) quantum system, otherwise known as a that a qu-$d$it (or qu$d$it)~\cite[2.2.1]{nielsenchuang:quantum}, with $\M_d(k)$ as its algebra of observables.
\item Given a quiver $\Gamma$, its \emph{path coalgebra} $k\Gamma$ is spanned by the paths in $\Gamma$, with comultiplication and counit defined on a path $p$ in $\Gamma$ of length $|p|$ by
\[
\Delta(p) = \sum_{p = p_1p_2} p_1 \otimes p_2 \quad \mbox{and} \quad
\epsilon(p) = \delta_{0,|p|}.
\]
If $\Gamma$ is a finite and acyclic, so that its path algebra $k[\Gamma]$ is finite-dimensional, then $k[\Gamma]^* \cong k\Gamma$. 
If $\Gamma$ has no arrows, the path coalgebra $k\Gamma$ coincides with the linearization of $\Gamma$ considered as a set. 
\item Let $\T_n = \bigoplus_{i\leq j} k E_{ij} \subseteq \M_n(k)$ denote the algebra of upper-triangular 
$n \times n$ matrices, and denote its dual coalgebra as $\T^n = \bigoplus_{i \leq j} E^{ij}$, so that 
we have a surjective coalgebra morphism $\M^n \twoheadrightarrow \T^n$. If we let $\Gamma$ be a quiver of type $A_n$ with linear orientation:
\[
\overset{1}{\bullet} \longleftarrow \overset{2}{\bullet} \longleftarrow \ \cdots \ \longleftarrow \overset{n}{\bullet}
\]
then $\T_n$ is isomorphic to the path algebra on $\Gamma$, so that its dual is isomorphic to the path coalgebra $\T^n \cong k\Gamma$. 
\end{enumerate}
\end{example}

The interpretation of coalgebras as quantized sets can also be seen at the level of dual algebras. 
The dual algebra of any coalgebra can be topologized by taking a neighborhood basis of zero as the vanishing ideals of any finite-dimensional subcoalgebra, resulting in a pseudocompact algebra~\cite[Section~3]{simson}.
If $X$ is a set, one can easily verify that the pseudocompact algebra dual to $kX$ is isomorphic to the product algebra
\begin{equation}\label{eq:functions}
(kX)^* \cong k^X
\end{equation}
(recalling that $k$ is discrete).
We may equivalently view $k^X$ as the algebra of $k$-valued functions on $X$,
in which case the topology on $k^X$ can be interpreted as the topology of pointwise convergence of functions.
Thus if we view a general coalgebra $Q$ as a noncommutative set, we can interpret $Q^*$ as the algebra of ``$k$-valued functions'' on $Q$, with its pseudocompact topology as a topology of ``quantum pointwise convergence.''

Viewing $Q^*$ as an algebra of functions on a quantum set $Q$ also suggests that there should be a correspondence between quantum discrete spaces (= quantum sets) and their function algebras, which is seen as follows. 
For ``classical'' sets, it is shown in~\cite[Theorem~4.7]{imr:diagonal} that the functor $X \mapsto k^X$ yields a duality between $\Set$ and a certain full subcategory $\catC$ of the category $\PCAlg_k$ of pseudocompact topological algebras. 
As described in~\cite[Theorem~3.6]{simson}, the dual algebra construction yields a dual equivalence
\[
(-)^* \colon \Coalg_k\op \to \PCAlg_k
\]
whose quasi-inverse is the continuous $k$-dual $(-)^\circ$.
In physical terms, we view this as a duality between state spaces and observables. The classical and quantum dualities thus fit into a diagram that commutes up to isomorphism
\[
\xymatrix{
\Set\, \ar[d]^[@!-90]{\sim} \ar@{^{(}->}[rr] & & \Coalg_k \ar[d]^[@!-90]{\sim}\\
\catC\op\, \ar@{^{(}->}[rr] & & \PCAlg_k\op
}
\] 
where the horizontal arrows are fully faithful embeddings.


Recall~\cite[Section~3.4]{radford:hopf} that the \emph{coradical} $\corad Q$ of a coalgebra $Q$ is defined to be the (direct) sum of all simple subcoalgebras, so that it is the largest cosemisimple subcoalgebra of $Q$. The span of the $k$-points of $Q$ forms a subcoalgebra of the coradical, so that we have coalgebra embeddings
\[
k \pts(Q) \hookrightarrow \corad Q \subseteq Q.
\]
If $k$ is algebraically closed, then the only cosemisimple coalgebras over $k$ are the comatrix coalgebras; in this case we have that $\corad Q$ consists of all qu$d$its in $Q$ while $k \pts(Q)$ is the span of all classical points. 
Thus if we wish to study coalgebras over more general fields, we can view the coradical as the collection of all ``generlized qu$d$its'' inside of a quantum set.
We will return to this interpretation in Subsection~\ref{sub:coradical} when discussing the coradical of a finite dual coalgebra.

\subsection{Coalgebras of distributions on commutative schemes}\label{sub:distributions}

To close this section, we explain how a certain coalgebra can be viewed as an ``underlying discrete object'' of
a scheme that is locally of finite type over $k$. 
This idea dates back to work of Takeuchi~\cites{takeuchi:tangent, takeuchi:formal}, where the \emph{underlying coalgebra} of a $k$-scheme is described in terms of representable functors.
By contrast, we construct these coalgebras in a more concrete manner via the language of \emph{distributions}, and we explicitly connect our presentation with that of Takeuchi in Proposition~\ref{prop:distribution sum}.
Note that a similar perspective in differential geometry is given by Batchelor in~\cite{batchelor:maps}.

We motivate our approach with the following observation about underlying sets of topological spaces.

\begin{remark}\label{rem:underlying set}
Let $X$ be a topological space satisfying the $T_1$ separation axiom, such as a Hausdorff space. Since points are closed in $X$, every finite subset $S \subseteq X$ is closed, and the subspace topology on $S$ is discrete. The family $\F(X)$ of finite subsets of $X$ is directed by inclusion. Taking the directed limit of this family in the  category of topological spaces yields a discrete space, which is the underlying set of $X$:
\[
\dirlim_{S \in \F(X)} S = |X|.
\]
\end{remark}

We will similarly consider ``underlying discrete objects'' of certain $k$-schemes $X$ by first
associating a suitable object to every closed subscheme of $X$ that is finite over $k$, and
then taking the direct limit of these objects. By contrast with the remark above, these ``discrete objects''
will be coalgebras that both ``linearize'' the set of closed points and include extra information
about the formal neighborhood of every closed point.

Let $\Sch_k$ denote the category of schemes over $k$, and let $\fSch_k$ denote the full subcategory of schemes that are finite over $\Spec(k)$, which we abbreviate in the typical way to \emph{finite over $k$}. 
Let $\cAlg$ denote the category of commutative $k$-algebras, and let $\fdcAlg$ denote the full subcategory of finite-dimensional commutative algebras. Similarly, let $\fdAlg$, $\fdCoalg$, and $\fdcCoalg$ respectively denote the categories of algebras, coalgebras, or cocommutative coalgebras that are finite-dimensional.

Let $S$ be a scheme finite over $k$, so that $S \cong \Spec(A)$ is affine, where $A$ is a 
finite-dimensional commutative $k$-algebra. 
We define the \emph{coalgebra of distributions} on $S$ to be the dual coalgebra
\[
\Dist(S) = \Gamma(S,\O_S)^*.
\]
As this is the composite of the global sections functor $\fSch_k\op \to \fdAlg$ with 
the contravariant dual functor $\fdAlg\op \to \fdCoalg$, this assignment yields a 
(covariant) functor
\begin{equation}\label{eq:finite distributions}
\Dist \colon \fSch_k \to \fdCoalg.
\end{equation}
In fact~\eqref{eq:finite distributions} is an equivalence of categories as it is the composition of the contravariant equivalences from $\fSch$ to finite-dimensional commutative algebras and from the latter to $\fdCoalg$ (see also~\cite[Theorem~1.1]{takeuchi:formal}).
%
%
In this sense, we imagine that finite-dimensional cocommutative coalgebras are ``the same as'' 
schemes finite over $k$. The covariance of the functor $\Dist$ suggests that we may view
these coalgebras as the ``underlying (discrete) object'' of such $k$-schemes.

We provide a basic example to which we will refer later. Consider the algebra $A = k[\eps]$ with $\eps^n = 0$ and basis $1,\eps, \dots, \eps^{n-1}$. Then $\Dist(\Spec(A)) \cong A^*$ is the coalgebra on the dual basis $\eps_r = \dual{\eps^r}$ having comultiplication and counit
\begin{equation}\label{eq:dual numbers}
\Delta\left( \eps_r \right) = \sum_{i + j = r} \eps_i \otimes \eps_j \quad \mbox{and} \quad
\epsilon\left( \eps_r \right) = \delta_{r0}.
\end{equation}
By the discussion above, we view this coalgebra as the ``underlying discrete object'' for the closed subscheme 
$\Spec(k[\eps]) \subseteq \A^1_k$. 

%

\begin{remark}\label{rem:finite image}
Let $g \colon S \to Y$ be a morphism in $\Sch_k$, and let $Z$ denote the scheme-theoretic image of $g$ in $Y$. If $S$ is finite over $k$, then the closed subscheme $Z$ of $Y$ is also finite over $k$. We omit the elementary proof of this fact, but make repeated use of it below.
\end{remark}

Now let $X$ denote an arbitrary $k$-scheme, and let $\F(X)$ denote the diagram of closed subschemes of $X$ that are finite over $k$ with the naturally induced closed immersions between them. We remark that $\F(X)$ is directed. 
Indeed, given $S \cong \Spec(B_1)$ and $T \cong \Spec(B_2)$ in $\F(X)$, one has the
union of closed subschemes $S \cup T \hookrightarrow X$ given by the ideal sheaf
$\mathcal{I}_{S \cup T} = \mathcal{I}_S \cap \mathcal{I}_T$. 
The immersions of $S$, $T$, and $S \cup T$ into $X$ factor through the coproduct $S \coprod T$ as
\[
\xymatrix@R=.1em{
S \ar[dr] & & & \\
& S \coprod T \ar@{->>}[r] & S \cup T \ar@{->>}[r] & X \\
T \ar[ur] & & &
}
\]
where $S \cup T$ is the scheme-theoretic image of the coproduct. 
Because $S \coprod T \cong \Spec(B_1 \oplus B_2)$ is finite over $k$, it follows
from Remark~\ref{rem:finite image} that $S \cup T$ is also finite over $k$. 

Applying the global sections functor to the diagram $\F(X)$ yields an inversely directed system 
of commuative algebras $\Gamma_S(S,\O_S)$ in $\fdAlg$, whose duals in turn form a 
directed system of cocommutative coalgebras in $\fdCoalg$. 
Motivated by Remark~\ref{rem:underlying set}, we arrive at the following definition of distributions for a general scheme.

\begin{definition}
For a $k$-scheme $X$, the \emph{coalgebra of distributions (of finite support)} on $X$ is the
directed limit of dual coalgebras
\[
\Dist(X) = \dirlim_{S \in \F(X)} \Dist(S) = \dirlim_{S \in \F(X)} \Gamma(S,\O_S)^*.
\]
\end{definition}

The assignment $X \mapsto \Dist(X)$ can be extended to a functor as follows. Let $f \colon X \to Y$
be a morphism of $k$-schemes. Given $S \in \F(X)$, we obtain a composite morphism of $k$-schemes
\[
S \hookrightarrow X \overset{f}{\to} Y.
\]
By Remark~\ref{rem:finite image}, the scheme-theoretic image $S'$ of $S$ in $Y$ is finite over $k$.
Thus (co)restriction of $f$ induces a morphism of $k$-schemes $S \to S'$, which induces a morphism 
of coalgebras $\Dist(S) \to \Dist(S')$ as in~\eqref{eq:finite distributions}. As $\Dist(X)$ is the 
colimit of the finite-dimensional subcoalgebras of the form $\Dist(S)$, and $\Dist(S') \hookrightarrow \Dist(Y)$
is a subcoalgebra, we may define
\[
\Dist(f) \colon \Dist(X) \to \Dist(Y)
\]
to be the directed limit of the morphisms induced from each $S \in \F(X)$.
In this way, distributions of finite support form a functor 
\begin{equation}\label{eq:distribution functor}
\Dist \colon \Sch_k \to \Coalg.
\end{equation}

In the case of an affine scheme over $k$, this functor amounts to the finite dual of the coordinate ring.

\begin{proposition}\label{prop:affine distributions}
The restriction of the functor $\Dist \colon \Sch_k \to \Coalg$ to the full subcategory of affine $k$-schemes is naturally isomorphic to the functor $X \mapsto \Gamma(X,\O_X)^\circ$. 
\end{proposition}

\begin{proof}
Writing $X \cong \Spec(A)$, it is clear that each $S \in \F(X)$ is of the form $S \cong \Spec(A/I)$ for some ideal $I \in \F(A)$ of finite codimension in $A$. Then 
\[
\Dist(X) = \dirlim_{S \in \F(X)} \Dist(S) \cong \dirlim_{I \in \F(A)} (A/I)^* = A^\circ \cong \Gamma(X, \O_X)^\circ.
\]
Naturality in $X$ is easily verified.
\end{proof}

A coalgebra of distributions supported at a point can also be defined in the following way. For a $k$-scheme $X$ and a point $x$ of $X$, the \emph{coalgebra of distributions
supported at $x$} is the finite dual of the stalk at $x$ of the structure sheaf:
\begin{equation}\label{eq:point distributions}
\Dist(X,x) = (\O_{X,x})^\circ.
\end{equation}
Note that if $x$ is not a closed point, then $\O_{X,x}$ has no ideals of finite codimension and $\Dist(X,x) = 0$; for this reason, distributions are typically only examined at closed points. 

Let $X_0$ denote the set of closed points of $X$. 
For $x \in X_0$, Takeuchi~\cite[2.1]{takeuchi:tangent} called~\eqref{eq:point distributions} the \emph{tangent coalgebra to $X$ at $x$} and referred to the direct sum $\bigoplus_{x \in X_0} \Dist(X,x)$ as the \emph{underlying coalgebra of $X$}. 
(See also~\cite[Chapter~7]{jantzen:groups} and~\cite[II, \S4, no.~5--6]{demazuregabriel} for details on coalgebras of  distributions at a point.)
Our next goal is to show that this underlying coalgebra is in fact isomorhpic to the coalgebra $\Dist(X)$ defined above.

\begin{example}\label{ex:distribution at 0}
We compute the coalgebra of distributions supported at the origin on the affine line over~$k$.
We have $\Dist(\A^1_k,0) = \dirlim (k[t]/(t^n))^*$, so that
\[
\Dist(\A^1,0) = \bigoplus_{i=0}^\infty k \eps_i
\] 
with the same comultiplication and counit formulas from~\eqref{eq:dual numbers}; this is known as the \emph{divided power coalgebra}.
(Note that this is isomorphic to the coalgebra associated to the monoid $\mathbb{N}$; see~\cite[p.~26]{radford:hopf}.)
Its dual algebra of ``functions'' is isomorphic to the formal power series ring 
$k[[t]] \cong \varprojlim k[t]/(t^n)$ with the $t$-adic topology. 
\end{example}

Let $x$ be a closed point of a $k$-scheme $X$. The canonical map $\Spec(\O_{X,x}) \to X$ induces a map on distributions
\[
\Dist(X,x) = (\O_{X,x})^\circ \cong \Dist(\Spec \O_{X,x}) \hookrightarrow \Dist(X),
\] 
and the assignment $(X,x) \mapsto \Dist(X,x)$ is evidently functorial (see also~\cite[Introduction]{takeuchi:tangent}).
In this way we obtain a naturally induced injection 
\begin{equation}\label{eq:dist sum}
\bigoplus_{x \in X_0} \Dist(X,x) \hookrightarrow \Dist(X).
\end{equation}
Note that the only nonzero summands above are from those points in the set
\[
\fin(X) = \{x \in X \mid [\kappa(x) : k] < \infty\} \subseteq X_0.
\]
If $X$ is locally of finite type over $k$, then we have $\fin(X) = X_0$ by the general Nullstellensatz.

We now verify our claim that the coalgebra of distributions on a $k$-scheme coincides with Takeuchi's underlying coalgebra. 

\begin{proposition}\label{prop:distribution sum}
For a $k$-scheme $X$, the map~\eqref{eq:dist sum} induces a natural isomorphism
of coalgebras 
\[
\Dist(X) \cong \bigoplus_{x \in \fin(X)} \Dist(X,x) = \bigoplus_{x \in \fin(X)} (\O_{X,x})^\circ.
\]
\end{proposition}

\begin{proof}
Let $S$ be a closed subscheme of $X$ that is finite over $k$. 
Let $\{x_1, \dots, x_n\}$ denote the underlying set of $S$. Consider each singleton $\{x_i\}$ as both a closed and an open subscheme of $S$ (since $S$ carries the discrete topology). This subscheme is isomorphic to $\Spec(B_i)$ for the finite-dimensional local algebra $B_i = \O_S(\{x_i\}) = \O_{S,x_i}$. Passing to an open affine neighborhood of $x_i$, we see that $B_i \cong \O_{X,x_i}/J_i$ for an ideal $J_i$ of finite codimension in the local ring $\O_{X,x_i}$. Notice that there exists a proper ideal of finite codimension in $\O_{X,x_i}$ if and only if $x_i \in \fin(X)$.

Setting $B = B_1 \times \cdots \times B_n$, we have the following commutative diagram
\[
\xymatrix{
\coprod \{x_i\} \ar[r]^-\sim \ar[d]^-[right]\sim & S \ar[d]^-[right]\sim \\
\coprod \Spec(B_i) \ar[r]^-\sim  & \Spec(B)
}
\]
of isomorphisms. On the level of distributions, we find that
\[
\Dist(S) \cong \bigoplus_{i=1}^n (B_i)^* \cong \bigoplus (\O_{X,x_i}/J_i)^*.
\]

If we now pass to the directed limit of all closed subschemes $S$ of $X$ that are finite over $k$, it is straightforward to see that
\[
\Dist(X) = \dirlim_{S \in \F(X)} \Dist(S) \cong \bigoplus_{x \in \fin(X)} (\O_{X,x})^\circ
\]
as desired.
\end{proof}

The closed points of a $k$-scheme are embodied in the structure of the coradical of $\Dist(X)$ in the following way.

\begin{corollary}\label{cor:distribution coradical}
Let $X$ be a $k$-scheme.
\begin{enumerate}
\item $\corad(\Dist(X)) = \bigoplus_{x \in X_0} \kappa(x)^\circ$, where $\kappa(x)$ denotes the residue
field at $x$.
\item \cite[2.1.8]{takeuchi:tangent} There is a bijection $\pts(\Dist(X)) \cong X(k)$ between the $k$-points of the distributions on $X$
and the $k$-rational points of $X$. 
\item If $k$ is algebraically closed and $X$ is locally of finite type over $k$, then $\corad(\Dist(X)) \cong k X_0$
(the coradical of $\Dist(X)$ is the linearization of the closed points of $X$). 
\end{enumerate}
\end{corollary}

\begin{proof}
Claim~(1) follows directly from the structure given in Proposition~\ref{prop:distribution sum}, and~(2) follows from~(1) because a closed point $x$ is $k$-rational iff $\kappa(x) = k$. Now~(3) follows from~(2) since every closed point is $k$-rational under the hypothesis.
\end{proof}

Assume that the scheme $X$ is locally of finite type over~$k$, so that $\fin(X) = X_0$. From Proposition~\ref{prop:distribution sum} and the corresponding isomorphism of pseudocompact algebras
\begin{equation}\label{eq:observables of distributions}
\Dist(X)^* \cong \prod_{x \in X_0} \widehat{\O}_{X,x},
\end{equation}
we see that $\Dist(X)$ contains information about the closed points of $X$ along with their formal neighborhoods.
While the coradical $\corad(\Dist(X))$ is a good choice of coalgebra corresponding to the closed points of $X$ (by Corollary~\ref{cor:distribution coradical}) and has good functorial properties for schemes, we will see in Theorem~\ref{thm:subfunctor} that it is necessary to retain the additional data of these formal neighborhoods in order for functoriality to persist in noncommutative geometry.

We conclude this section by describing the distributions on the affine line over an algebraically closed field. This coalgebra will appear again when examining the quantum plane in Subection~\ref{sub:quantum plane}.

\begin{example}\label{ex:affine line}
Suppose $k$ is algebraically closed, so that the closed points of $\A^1_k = \Spec k[t]$ are in bijection with $k$ and are all $k$-rational. From Proposition~\ref{prop:distribution sum} we obtain 
\begin{equation}\label{eq:line decomposition}
\Dist(\A^1_k) \cong \bigoplus_{\lambda \in k} \Dist(\A^1_k,\lambda).
\end{equation}
By translation invariance we have each $\Dist(\A^1_k, \lambda) \cong \Dist(\A^1, 0)$, which was already described in Example~\ref{ex:distribution at 0}. However, it will be important later for us to describe this coalgebra with a more careful consideration of the pointwise structure. At each point $\lambda \in k$, we have 
\begin{equation}\label{eq:}
\Dist(\A^1_k, \lambda) = \dirlim_n (k[t]/(t-\lambda)^n)^* =  \bigoplus_{i=0}^\infty k\eps^{(i)}_\lambda,
\end{equation}\label{eq:line at a point}
and its coalgebra structure is given by 
\begin{align*}
\Delta(\eps_\lambda^{(r)}) &= \sum_{i=0}^r \eps_\lambda^{(i)} \otimes \eps_\lambda^{(r-i)}, \\
\epsilon(\eps_\lambda^{(r)}) &= \delta_{r0}.
\end{align*}
Viewing these $\eps_\lambda^{(i)}$ as distributions using~\eqref{eq:line decomposition} and Proposition~\ref{prop:affine distributions}, one can check that $\eps_\lambda^{(i)} \in k[t]^\circ$ is described on the basis $\{(t-\lambda)^j\}_{j = 0}^\infty$ by $\eps_\lambda^{(i)}((t-\lambda)^j) = \delta_{ij}$.
In order to find a describe $\Dist(\A^1_k)$ in terms of distributions that are independent of a choice of basis of $k[t]$, it is instructive to rewrite this structure in terms of Dirac distributions~\cite[Example~2.1.2]{hormander} in the characteristic zero case.
The Dirac distribution $\delta_\lambda \in \Dist(\A^1_k)$ supported at $\lambda \in k$ is the functional that evaluates at $t = \lambda$:
\[
\delta_\lambda(f(t)) = f(\lambda).
\]
In particular, $\delta_\lambda = \eps_\lambda^{(0)}$.
The distributional derivatives~\cite[Definition~3.1.1]{hormander} of the Dirac distributions $\delta = \delta_\lambda$ are defined by
\begin{align*}
\delta'(f) &= - \delta(f'), \mbox{ and more generally}\\
\delta^{(i)}(f) &= (-1)^i \delta(f^{(i)}).
\end{align*}
This means that if $k$ has characteristic~0,
\[
\eps^{(i)}_\lambda = \frac{(-1)^i}{i!} \, \delta_\lambda^{(i)}.
\]
So  the comultiplication $\Delta(\eps_\lambda^{(r)}) = \sum \epsilon_\lambda^{(i)} \otimes \epsilon_\lambda^{(r-i)}$ translates to
\begin{align*}
\Delta(\delta_\lambda^{(r)}) 
&= \frac{r!}{(-1)^r} \sum_{i=0}^r \frac{(-1)^i}{i!} \delta_\lambda^{(i)} \otimes \frac{(-1)^{r-i}}{(r-i)!}\delta_\lambda^{(j)} \\
&=  \sum_{i=0}^r \binom{r}{i} \delta_\lambda^{(i)} \otimes \delta_\lambda^{(r-i)}.
\end{align*}
\end{example}

\section{The finite dual as a quantized maximal spectrum}\label{sec:quantized Max}

Let $A$ be a $k$-algebra. If $A$ happens to be a commutative affine $k$-algebra, then the corresponding scheme $X = \Spec A$ is of finite type over~$k$, and
Propositions~\ref{prop:affine distributions} and~\ref{prop:distribution sum} combine to give the following isomorphisms that are natural in $A$:
\begin{equation}\label{eq:quantum Max}
\Dist(\Spec(A)) \cong A^\circ \cong \bigoplus_{\m \in \Max(A)} (A_\m)^\circ.
\end{equation}
(The pseudocompact dual of this isomorphism was noted in~\cite[Example~2.9]{kontsevichsoibelman} for algebraically closed fields, and more generally in~\cite[Section~2]{lebruyn:dual}.)
This isomorphism leads directly to the main tenet of this paper: 
\begin{quote}
\emph{For affine $k$-algebras that are not necessarily commutative but have enough finite-dimensional $k$-representations, the finite dual coalgebra functor $(-)^\circ$ is a meaningful substitute for the maximal spectrum functor.} 
\end{quote}

Another way to motivate this is by viewing measuring coalgebras $P(A,B)$ as a quantization~\cite{batchelor:ncg} of the set of maps between two $k$-algebras $A$ and $B$. This provides an enrichment of the category $\Alg$ of $k$-algebras over the category $\Coalg$; see~\cite{vasilakopoulu:enrichment, hlfv:enrichment}. If $k$ is algebraically closed and $A$ is a commutative affine algebra, then the Nullstellensatz yields a bijection
\[
\Alg(A,k) \cong \Max(A).
\]
If we upgrade the above Hom-set using this enriched Hom structure, it is known~\cite[Remark~3.1]{vasilakopoulu:enrichment} that
\[
P(A,k) \cong A^\circ.
\]
This provides a second perspective on the finite dual as a quantized maximal spectrum. Note that the adjoint relationship~\eqref{eq:adjoint} allows us to recover the former Hom-set from the enriched one using $k$-points, as
\[
\pts(A^\circ) = \Coalg(k,A^\circ) \cong \Alg(A,k).
\]

The dual algebra of ``functions'' on a dual coalgebra can be described as follows. Since $A^\circ = \dirlim_{I \in \F(A)} (A/I)^*$ with notation as in subsection~\ref{sub:Sweedler}, we have
\[
(A^\circ)^* \cong \invlim_{I \in \F(A)} A/I =: \widehat{A}.
\]
In the case where $A$ is commutative and affine, then as in~\eqref{eq:observables of distributions} we have $\widehat{A} \cong \prod_{\m \in \Max(A)} \widehat{A_\m}$. This underscores the fact that $A^\circ$ contains information not only about closed points, but also about their infinitesimal neighborhoods. We will see in Theorem~\ref{thm:subfunctor} that this infinitesimal information is unavoidable if we desire a useful noncommutative spectrum functor.

The remainder of this section is divided as follows. First we will discuss the failure of functoriality of the coradical of $A^\circ$ to clarify our insistence on retaining the full finite dual when quantizing the maximal spectrum functor; we also examine a few cases where the coradical happens to be functorial. We then discuss those noncommutative algebras $A$ for which we believe $A^\circ$ is an appropriate quantization of the maximal spectrum, the \emph{fully residually finite-dimensional} algebras. We end by examining the relationship between Morita equivalence of algebras and Takeuchi equivalence of their dual coalgebras, treating it as an object lesson on the role of Morita equivalence in functorial spectral theory.

\subsection{Functoriality and the coradical}\label{sub:coradical}

In light of Corollary~\ref{cor:distribution coradical}, it seems natural to infer from the isomorphism~\ref{eq:quantum Max} that the coradical of $A^\circ$ is a more appropriate quantization of the maximal spectrum of a commutative affine algebra. 
It is well known (see also Corollary~\ref{cor:functorial cases}(1) below) that the coradical of the finite dual forms a functor on commutative $k$-algebras:
\begin{align}\label{eq:commutative coradical}
\cAlg\op &\to \Coalg \\
A &\mapsto \corad A^\circ. \nonumber
\end{align}
Indeed, we agree that for suitable algebras $A$ (as discussed in Subsection~\ref{sub:fully RFD} below) the coradical of $A^\circ$ is an appropriate substitute for the set of points of the noncommutative affine variety associated to $A$. What we illustrate in this subsection is rather that, in order to preserve functoriality of the maximal spectrum in these cases, one is forced to include the full finite dual. So while points of affine varieties may behave well in commutative algebraic geometry, we are forced to include infinitesimal neighborhoods in order to retain functoriality in noncommutative geometry.

First let us describe the structure of $\corad A^\circ$. It is clear from the isomorphism~\eqref{eq:Sweedler colimit} that the simple subcoalgebras of $A^\circ$ correspond to the quotient algebras $A/I$ that are simple and finite-dimensional. Let
\[
\{\m_\alpha\} = \Max(A) \cap \F(A)
\]
be an indexing of the set of maximal ideals of finite codimension in $A$, and let $S_\alpha = A/\m_\alpha$ denote the simple quotient algebras. Then we have
\begin{equation}\label{eq:qubit part}
\corad A^\circ = \bigoplus_{\alpha} S_\alpha^*.
\end{equation}
If $S_\alpha \cong \M_d(k)$, then as discussed in Subsection~\ref{sub:quantized sets} we may view $S_\alpha^* \cong \M^d$ as a qu$d$it over~$k$. So $S_\alpha^*$ might generally be imagined as representing a ``generalized qu$d$it'' over~$k$ (although this is probably most appropriate in the case where $S_\alpha$ is central simple over $k$).

Thus in light of~\eqref{eq:qubit part}, we may interpret $\corad A^\circ$ as a ``disjoint union of generalized qu$d$its'' in the noncommutative spectrum of $A$. (Note that if $k$ is algebraically closed, then all $S_\alpha$ are isomorphic to matrix algebras, so that the coradical truly does represent a disjoint union of qu$d$its over $k$.)
Furthermore, there is a bijection between the $\m_\alpha$ and the finite-dimensional simple left $A$-modules ${}_A V_\alpha$, where $\m_\alpha$ is the annihilator of $V_\alpha$. So $\corad A^\circ$ essentially encodes the irreducible finite-dimensional representations of $A$.

We now turn to the issue of non-functoriality of the coradical.
To begin, very simple examples can be constructed to show that the coradical does not form a subfunctor of the finite dual, even when restricting to the full subcategory of affine noetherian PI algebras that are module-finite over their centers.
We find the following particularly helpful as an instance that can be visualized.

\begin{example}
Consider the algebra homomorphism $\phi \colon k[t] \to \M_2(k)$ given by $\phi(t) = E_{12}$. Its image is isomorphic to the ring of dual numbers $k[\eps]$ (with $\eps^2 = 0$), so that $\phi$ factors as $k[t] \twoheadrightarrow k[\eps] \hookrightarrow \M_2(k)$. This dualizes to a factorization of coalgebra morphisms $\M^2 = \M_2(k)^* \twoheadrightarrow k[\eps]^* \hookrightarrow k[t]^\circ$. In particular, we see that the image of $\corad(\M^2) = \M^2$ under $\phi^\circ$ is not cosemisimple and thus does not lie in $\corad k[t]^\circ$.

To paint a more suggestive picture, let us first note that the image of $\phi$ in lies in the subalgebra $\T_2(k) \subseteq \M_2(k)$ of upper-triangular matrices, so that $\phi$ factors as
\[
k[t] \twoheadrightarrow k[\eps] \hookrightarrow \T_2(k) \hookrightarrow \M_2(k). 
\]
The finite duals of the algebras were described in Examples~\ref{ex:coalgebra examples} and~\ref{ex:affine line}. 
So the above factorization of $\phi$ dualizes to a factorization of coalgebra morphisms
\begin{equation}\label{eq:qubit map}
\M^2 \twoheadrightarrow \T^2 \twoheadrightarrow k[\eps]^\circ \hookrightarrow  \Dist(\A^1_k).
\end{equation}
Recall from Example~\ref{ex:coalgebra examples} that if we let $\Gamma$ be a quiver of type $A_2$, then we may view $\T^2 = k\Gamma$ as a path coalgebra.
Figure~\ref{fig:qubit map} is a visualization the sequence of morphisms above. We use the traditional Bloch sphere depiction~\cite[p.~15]{nielsenchuang:quantum} of the qubit to represent $\M^2$, the quiver $\Gamma$ to depict $\T^2 = k\Gamma$, and an infinitesimal neighborhood of the origin to represent $k[\eps]^*$.
\begin{figure}
\includegraphics[scale=.7]{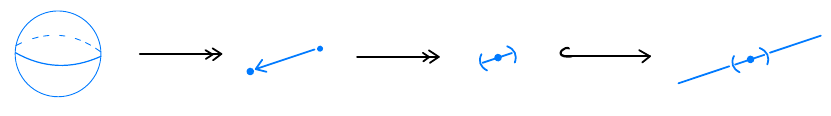}
\caption{Illustration of the maps between quantum sets in~\eqref{eq:qubit map}}
\label{fig:qubit map}
\end{figure}
In total, we have a map from the qubit to the affine line (both defined over~$k$) whose image is a non-reduced subscheme of $\A^1_k$. 
\end{example}

In the spirit of functorial noncommutative spectral theory as described in Subsection~\ref{sub:spectral}, we prefer stronger evidence than a single example in order to conclude that the full finite dual functor should not be replaced by a functor 
\[
F \colon \Alg\op \to \Coalg
\]
that is a better approximation to the coradical. Theorem~\ref{thm:subfunctor} below gives a more precise restriction.

To motivate the statement, we first note that there exist functors $F$ as above whose restriction to $\cAlg\op$ is the functor~\eqref{eq:commutative coradical}.
Given a $k$-algebra $A$, let $F(A) \subseteq \corad A^\circ$ be the sum of all simple subcoalgebras dual to a a finite-dimensional division algebra, which is the largest basic~\cite[p.~43]{chinmontgomery:basic} subcoalgebra. The argument given in the proof of Corollary~\ref{cor:functorial cases}(1) below shows that  this choice of $F$ is a subfunctor of the finite dual, and it is routine to check that 
$F(A) = \corad A^\circ$
for commutative algebras $A$. 
However, this functor is easily seen to satisfy $F(\M_d(k)) = 0$ for $d \geq 2$, so that matrix algebras yield familiar obstructions as in~\cites{reyes:obstructing, bergheunen:obstructions, reyes:sheaves, bmr:kochenspecker}. 

So from the point of view of Subsection~\ref{sub:spectral}, we wish to impose a nondegeneracy condition which guarantees that, at the very least, matrix algebras are assigned nontrivial coalgebras.
If we resolve to view matrix algebras as observables on qu$d$its as discussed in Subsection~\ref{sub:quantized sets}, then it is more natural to seek a subfunctor $F$ of $(-)^\circ$ that satisfies $F(\M_d(k)) =  \M_d(k)^* = \M^d$ for all matrix algebras. This leads us uniquely to the functor $F = (-)^\circ$ in the following way.

\begin{theorem}\label{thm:subfunctor}
Let $F \colon \Alg\op \to \Coalg$ be a subfunctor of the finite dual~\eqref{eq:Sweedler}. If $F(\M_d(k)) = \M_d(k)^*$ for all integers $n \geq 0$, then $F$ is equal to the finite dual.
\end{theorem}

\begin{proof}
First let $A$ be a finite-dimensional $k$-algebra, and set $d = \dim_k(A)$. Consider the embedding of $A$ into its $k$-endomorphisms:
\[
i \colon A \hookrightarrow \End_k(A) =: S.
\]
Because $S \cong \M_d(k)$, our hypothesis yields $F(S) = S^*$. As $F$ is a subfunctor of $(-)^\circ$, we obtain a commuting diagram
\[
\xymatrix{
F(S)\ar@{=}[d] \ar[r]^{F(i)} & F(A) \ar@{^{(}->}[d]^{\alpha} \\
S^* \ar@{->>}[r]^{i^*}& A^*
}
\]
The morphism $\alpha$ is an inclusion of a subcoalgebra by the subfunctor condition. Now a simple diagram chase incorporating surjectivity of $i^*$ shows that $\alpha$ is also surjective, and therefore $F(A) = A^*$ (and $\alpha$ is the identity). 

Since $F$ is a subfunctor of the finite dual, we now deduce that $F = (-)^*$ on the full subcategory of finite-dimensional $k$-algebras. It is then straightforward to apply the directed colimit characterization~\eqref{eq:Sweedler colimit} of the finite dual to deduce that in fact $F = (-)^\circ$ on the category $\Alg$ of all $k$-algebras.
\end{proof}

This is one way in which we might think of the finite dual as a ``minimal'' functorial extension of the coradical from commutative to noncommutative algebras. Let us return to the the description~\eqref{eq:qubit part} of the coradical of $A^\circ$ once more. As discussed above, the coradical contains information about all finite-dimensional irreducible representations of $A$. While $\corad A^\circ$ fails to behave functorially in general, the theorem above suggests that we can view the finite dual as a minimal functorial substitute for the irreps of $A$. Or stated otherwise, $A^\circ$ provides us with a ``functorial snapshot'' of the representation theory of $A$.

Despite the general failure of functoriality for the coradical, there are conditions on an algebra homomorphism under which this coradical \emph{does} behave functorially. These can be characterized in the following way.

\begin{proposition}\label{prop:functorial coradical}
For an algebra homomorphism $f \colon A \to B$, the following are equivalent:
\begin{enumerate}[label=\textnormal{(\alph*)}]
\item $f^\circ(\corad B^\circ) \subseteq \corad A^\circ$ (i.e., $f^\circ$ preserves the coradical);
\item For every finite-dimensional semisimple left (resp., right) $B$-module $M$, the restriction of scalars ${}_A M$ (via $f$) is a semisimple $A$-module;
\item For every maximal ideal $\m$ of finite codimension in $B$, the ideal $f^{-1}(\m) \subseteq A$ is semiprime. 
\end{enumerate}
\end{proposition}

\begin{proof}
The simple subcoalgebras of $B^\circ$ are precisely those of the form $(B/\m)^*$ for a maximal ideal $\m$ of finite codimension in $B$. The restriction of $f^\circ$ to one of these simple subcoalgebras coincides with the surjective map 
\[
\Hom_k(B/\m,k) \to \Hom_k(A/f^{-1}(\m),k)
\]
arising from the algebra homomorphism $A/f^{-1}(\m) \to B/\m$ induced by $f$. Thus the image of this simple subcoalgebra of $B^\circ$ lies in the coradical of $A^\circ$ if and only if the finite-dimensional algebra $A/f^{-1}(\m)$ is semisimple, which occurs if and only if $f^{-1}(\m)$ is semiprime. This establishes (a)$\iff$(c).

To see that (b)$\implies$(c), let $\m$ be a maximal ideal of finite codimension in $B$. Then $B/\m$ is a semisimple left $B$-module, and under hypthesis~(b) it is also a semisimple left $A$-module. But the natural ring embedding $A/f^{-1}(\m) \hookrightarrow B/\m$ is also an embedding of left $A$-modules. So $A/f^{-1}(\m)$ is semisimple and consequently $f^{-1}(\m)$ is semiprime.

Finally, assume~(c). For (c)$\implies$(b), it suffices to consider the case where ${}_B M$ is simple and prove that ${}_A M$ is semisimple. Setting $\m = \ann_B(M)$, it follows that $B/\m$ is a finite-dimensional simple algebra. Because $\ann_A(M) = f^{-1}(\m)$, the action of $A$ on $M$ factors through $A/f^{-1}(\m)$. This finite-dimensional algebra is semiprime by hypothesis, so it must be semisimple. Thus $M$ is semisimple as a left $A$-module.
\end{proof}

We remark that if $f \colon A \to B$ satisfies the equivalent conditions above, then it follows from~\cite[Corollary~4.2.2]{radford:hopf} that if $f^\circ$ preserves the whole coradical filtration.
Below are two situations in which the algebra $B$ is ``tame enough'' for this to automatically hold.

\begin{corollary}\label{cor:functorial cases}
An algebra homomorphism $f \colon A \to B$ satisfies $f^\circ(\corad B^\circ) \subseteq \corad A^\circ$ in each of the following cases:
\begin{enumerate}
\item Every finite-dimensional simple quotient algebra of $B$ is a division algebra. (This holds, for instance, if $B$ is commutative.)
\item $B$ is a finite normalizing extension~\cite[10.1.3]{mcconnellrobson} of the image $f(A)$. 
\end{enumerate}
\end{corollary}

\begin{proof}
(1) This case is a very slight generalization of~\cite[Exercise~4.1.2]{radford:hopf}. If $\m$ is a maximal ideal of finite codimension $B$, then the image of $f(A)$ in $B/\m$ is a finite-dimensional domain and therefore is a division algebra. So $f^{-1}(\m)$ is maximal in $A$, and the result follows from Proposition~\ref{prop:functorial coradical}(c).

(2) This case also follows from Proposition~\ref{prop:functorial coradical}(c), this time with the help of the ``Cutting Down'' theorem~\cite[Theorem~10.2.4]{mcconnellrobson}.
\end{proof}

Note that condition~(1) is satisfied when $B$ is commutative, or more generally, when $B/J(B)$ is commutative where $J(B)$ is the Jacobson radical. Condition~(2) is satisfied in case $B$ is module-finite over a central subalgebra $A$ and $f$ is the inclusion of $A$ into $B$. An important case of the latter occurs when $B$ is a Cayley-Hamilton algebra and $A$ is the image of its trace; see~\cite[Theorem~4.5]{deconciniprocesi:quantum}.

\subsection{Fully RFD algebras}\label{sub:fully RFD}

We now discuss the question of which algebras have ``enough'' finite-dimensional representations that $A^\circ$ can be considered as a reasonable substitute for a noncommutative maximal spectrum. If the algebra structure of $A$ is to be faithfully represented by $A^\circ$, we should at least ask that the naturally induced algebra homomorphism
$A \to A^{\circ *}$
is injective. Injectivity of this map is equivalent~\cite[Lemma~6.1.0]{sweedler:hopf} to the property that the algebra is a subdirect product of finite-dimensional algebras. A $k$-algebra satisfying these conditions is said to be \emph{residually finite-dimensional} (RFD). It is straightforward to verify that the RFD property is further equivalent to the requirement that intersection of all ideals in $\F(A)$ is zero, or that every nonzero element of $A$ acts nontrivially on a finite-dimensional left (equivalently, right) module. 

However, it seems prudent to ask for more than the property above. For instance, if $I$ is an ideal of $A$ and if we view $A/I$ as the coordinate ring of a closed subscheme of the spectrum of $A$, then we should reasonably expect $A^\circ$ to carry sufficient information about $A/I$ as well. This will not always be the case, as RFD algebras can have homomorphic images that are not RFD. For instance, the algebra $A = \prod_{n=1}^\infty \M_n(k)$ is RFD by construction, but for the ideal $T = \bigoplus_{n=1}^\infty \M_n(k)$ of $A$, the same argument as in the proof of~\cite[Lemma~7.5]{kaplansky:operator} shows that $A/T$ has no finite-dimensional homomorphic image. Similarly, the free algebra $k\gen{x,y}$ is known to be RFD (as the ideals $I_n$ generated by words in $\{x,y\}$ of length~$n$ have finite codimension and satisfy $\bigcap I_n = 0$), but its homomorphic image $k\gen{x,y}/(xy-1)$ is not Dedekind-finite and therefore cannot be a subdirect product of finite-dimensional algebras.

\begin{definition}
(Following~\cite{courtneyshulman} in the case of C*-algebras.) An algebra $A$ is \emph{strongly residually finite-dimensional} (strongly RFD) if $A/I$ is residually finite-dimensional for every ideal $I$ of $A$. 
\end{definition}

Every affine noetherian PI $k$-algebra $A$ is strongly RFD. Indeed, a result of Anan'in~\cite{ananin:representability} shows that every affine right noetherian PI algebra is RFD. Because the property of being an affine noetherian PI algebra passes to homomorphic images, it follows that $A/I$ is RFD for all ideals $I$ of $A$.

From the discussion above, we see that $A$ is strongly RFD if and only if every ideal $I$ of $A$ is equal to the intersection of all ideals of finite codimension containing $I$, if and only if every homomorphic image of $A$ embeds in the algebra of observables on its finite dual. 
Thus the algebras for which $A^\circ$ can be considered a reasonable quantization of the maximal spectrum should form a subclass of the strongly RFD algebras. However, one could reasonably ask for a slightly stronger condition. Because $A^\circ$ contains information about finite-dimensional representations of $A$ as in~\eqref{eq:qubit part}, we can ask that every finitely generated $A$-module be determined by all of its finite-dimensional homomorphic images. This reasoning leads to the following condition.

\begin{definition}
A $k$-algebra $A$ is \emph{left fully residually finite-dimensional} (RFD) if every finitely generated left $A$-module is a subdirect product of finite-dimensional left $A$-modules. \emph{Right fully residually finite-dimensional} algebras are defined similarly, and we say that $A$ is \emph{fully residually finite-dimensional} if it is both left and right fully RFD.
\end{definition}

In general, we have the following relationship between these properties of $k$-algebras:
\[
\mbox{left fully RFD} \implies \mbox{strongly RFD} \implies \mbox{RFD}.
\]
We have seen that the second implication is strict. We presume that the first implication is also strict in general, but we do not have a counterexample. Proposition~\ref{prop:FBN RFD} below will illustrate that any such counterexample cannot be a fully bounded noetherian algebra.

Below are some equivalent characterizations of the left fully RFD property. Condition~(c), in particular, highlights the role of the finite dual for such algebras.
A left $A$-module is \emph{locally finite} if each of its finitely generated submodules is finite-dimensional. We let $A \lMod_{lf}$ denote the full subcategory of $A \lMod$ whose objects are the locally finite modules.

\begin{theorem}\label{thm:fully RFD}
For an algebra $A$, the following are equivalent:
\begin{enumerate}[label=\textnormal{(\alph*)}]
\item $A$ is left fully RFD; 
\item the injective hull of every simple left $A$-module is locally finite-dimensional;
\item $A^\circ$ is a cogenerator in $A \lMod$.
\end{enumerate}
If $A$ is left noetherian, then these properties are further equivalent to:
\begin{enumerate}[resume, label=\textnormal{(\alph*)}]
\item every simple left $A$-module is finite-dimensional, and the subcategory $A \lMod_{lf}$ is closed under injective hulls in $A \lMod$;
\item $A^\circ$ is an injective cogenerator in $A \lMod$;
\end{enumerate}
\end{theorem}

\begin{proof} 
Let $\{V_i\}$ be a complete set of simple left $A$-modules up to isomorphism.
Because every left $A$-module is a subdirect product of the injective hulls of the simple left $A$-modules (\cite[Theorem~19.8]{lam:lectures}), every finitely generated left $A$-module is a subdirect product of finitely generated submodules of the $E(V_i)$. Furthermore, each $E(V_i)$ and its submodules are subdirectly irreducible (since they contain an essential simple submodule). The equivalence (a)$\iff$(b) follows.

(b)$\implies$(c): Let $V$ be a simple left $A$-module, and let $\m = \ann(V)$. Because the injective hull $E(V)$ is locally finite-dimensional, $V$ must be finite-dimensional. Thus the semisimple algebra $\End_A(V) \cong A/\m$ is also finite-dimensional and consequently is symmetric~\cite[Example~16.59]{lam:lectures}. It follows that we have an embedding of finite-dimensional left $A$-modules
\[
V \hookrightarrow A/\m \cong (A/\m)^* \hookrightarrow A^\circ.
\]
As explained in the proof of~\cite[Theorem~2.1]{hatipoglulomp:hopf}, the left module $A^\circ$ is an injective object in the category locally finite left $A$-modules. Since~(b) states that $E(V)$ is locally finite-dimensional, the embedding above must extend to an injective homomorphism $E(V) \hookrightarrow A^\circ$. It follwos~\cite[Theorem~19.8]{lam:lectures} that $A^\circ$ is a cogenerator in $A \lMod$. 

(c)$\implies$(b): Let $V$ be a simple left $A$-module. If $A^\circ$ is a cogenerator in $A \lMod$, then there is an injective homomorphism $E(V) \hookrightarrow A^\circ$ (see~\cite[Theorem~19.8]{lam:lectures}). Because $A^\circ$ is locally finite-dimensional, the same must be true for $E(V)$.

(d)$\implies$(e): As explained in the proof of~\cite[Theorem~2.1]{hatipoglulomp:hopf}, if $A \lMod_{lf}$ is closed under injective hulls then $A^\circ$ is an injective left $A$-module (even without the left noetherian hypothesis). It is also clear that (d)$\implies$(b), so from (b)$\implies$(c) established above we have that $A^\circ$ is also a cogenerator.

Obviously (e)$\implies$(c). 
Now assume $A$ is left noetherian; we will verify (b)$\implies$(d). Suppose $M$ is a locally finite left $A$-module. Let $\{M_i\}$ denote the set of all finite-dimensional submodules of $M$. Then $M = \dirlim M_i$ is the directed union of the $M_i$. Because $A$ is left noetherian, injective left $A$-modules are closed under direct limits~\cite[Theorem~3.46]{lam:lectures} so that
\[
E(M) = \dirlim E(M_i).
\]
Since each $M_i$ has essential socle, its injective hull $E(M_i)$ is a finite direct sum of injective hulls of simple modules, which are locally finite-dimensional. Thus $E(M)$ is a directed union of locally finite-dimensional modules and must itself be locally finite-dimensional.
\end{proof}

This raises the obvious question of how to locate examples of fully RFD algebras. Fortunately, it turns out that there is a rich supply of these, in the sense that \emph{every affine noetherian PI algebra is fully RFD.} This is deduced in the next result, which also shows that for fully bounded noetherian (FBN) rings~\cite[Chapter~9]{goodearlwarfield}, the left and right fully RFD properties are equivalent to one another. 

\begin{proposition}\label{prop:FBN RFD}
Let $A$ be a fully bounded noetherian $k$-algebra. The following are equivalent:
\begin{enumerate}[label=\textnormal{(\alph*)}]
\item $A$ is fully RFD;
\item $A$ is strongly RFD;
\item every simple left (equivalently, right) $A$-module is finite-dimensional;
\item every maximal ideal of $A$ has finite codimension.
\end{enumerate}
In particular, every affine noetherian PI algebra is fully RFD. 
\end{proposition}

\begin{proof}
The implications (a)$\implies$(b)$\implies$(c) hold for any algebra $A$.
We use the following results on FBN rings:
\begin{enumerate}[label=(\roman*)]
\item The annihilator of every simple left or right module is maximal~\cite[Corollary~9.5]{goodearlwarfield}.
\item The injective hull of a simple left or right module is locally of finite length~\cite[Theorem~3.5, Corollary~3.6]{jategaonkar:jacobsons}.
\end{enumerate}
The equivalence (c)$\iff$(d) follows directly from~(i). Thanks to~(ii) and Theorem~\ref{thm:fully RFD}(c), $A$ is left fully RFD if and only if every simple left $A$-module has finite dimension. Since the left and right module versions of~(c) are equivalent, it follows that (a)$\iff$(c). 

The final statement is a consequence of the following two facts: every PI~algebra is fully bounded~\cite[Corollary~13.6.6]{mcconnellrobson}, and condition~(c) holds over any affine PI algebra~\cite[Theorem~13.10.3]{mcconnellrobson}.
\end{proof}

\begin{remark}
By a result of Amitsur and Small~\cite{amitsursmall}, every affine FBN algebra over an uncountable algebraically closed field~$k$ is PI, from which it follows that such a $k$-algebra is fully RFD. It seems natural to then ask if every affine FBN algebra over an arbitrary field~$k$ is fully RFD over~$k$. However, because division algebras are FBN, this question is in fact a generalization of a well-known open problem~\cite[Question~5]{smoktunowicz} which asks if every affine division algebra is finite-dimensional. 
Indeed, if $A$ is an affine FBN algebra with maximal ideal $M$, then $A/M \cong \M_n(D)$ for a division algebra $D$ that must be affine by the Artin-Tate Lemma. Assuming this problem has a positive solution, $D$ is finite-dimensional so that $M$ has finite codimension and $A$ is fully RFD by Proposition~\ref{prop:FBN RFD}.
\end{remark}

There are many more interesting questions that one can ask about the fully residually finite-dimensional property. 
Are the left and right fully RFD properties independent of one another? If $A^\circ$ is an \emph{injective} cogenerator in $A \lMod$, must $A$ be left noetherian? Are there any examples of affine but non-noetherian fully RFD algebras? How similar is the structure theory of fully RFD algebras to that of FBN rings? For now, we must content ourselves with the information learned above and turn to other matters regarding the finite dual.

\subsection{Morita equivalence and the quantized spectrum}

A common philosophy within noncommutative geometry is that \emph{Morita equivalent rings should represent the same noncommutative space.} The claim is natural and understandable when one works in a framework where a (commutative or noncommutative) space is represented solely by its category of sheaves of modules. However, this is problematic from the perspective laid out in Section~\ref{sec:intro}. For instance, if we wish to view a noncommutative algebra as consisting of observables of a quantum system, then its commutative subalgebras are physically important as they are closely related to the information that can be accessed through all possible measurements~\cite[1.3]{heunen:classicalfaces}. In particular, the algebra $\M_n(\C)$ of observables on an $n$-level system has many commutative subalgebras isomorphic to $\C^n$, representing different measurements with $n$ possible outcomes that can be made on the system. However, all of these matrix algebras $\M_n(\C)$ are Morita equivalent to one another and, in particular, to $\C$. Thus Morita equivalence is blind to this important aspect of algebraic quantum mechanics, while functorial spectral theory as in Subsection~\ref{sub:spectral} aims to preserve exactly this kind of information.  

This is not to say that Morita equivalence should be ignored in functorial spectral theory. Rather, from this perspective we expect that \emph{Morita equivalent algebras will have Morita equivalent noncommutative spaces, but Morita equivalent noncommutative spaces need not be isomorphic.} This subsection illustrates this principle in action for dual coalgebras. Two coalgebras are said to be \emph{Takeuchi equivalent}~\cite{takeuchi:morita} if their categories of left comodules are $k$-linearly equivalent. We will show below that if $A$ and $B$ are algebras that are Morita equivalent (in an appropriately modified sense), then $A^\circ$ and $B^\circ$ are Takeuchi equivalent. So this quantized maximal spectrum reflects Morita equivalence, in accordance with the guideline above.

We will say that two $k$-algebras $A$ and $B$ are \emph{$k$-linearly Morita equivalent} if there is a $k$-linear equivalence between the categories $A \lMod$ and $B \lMod$. A number of standard facts~\cite[Section~18]{lam:lectures} about Morita equivalence of rings are readily verified to carry over to $k$-linear Morita equivalence as long as suitable care is given to compatibility with the $k$-vector space structure. In particular, one can verify that within the class of $k$-algebras, the RFD, strongly RFD, and (left) fully RFD properties
are each preserved by $k$-linear Morita equivalence.

\begin{proposition}
Let $A$ and $B$ be $k$-algebras, and consider the following statements:
\begin{enumerate}[label=\textnormal{(\alph*)}]
\item $A$ and $B$ are $k$-linearly Morita equivalent algebras;
\item $A^\circ$ and $B^\circ$ are Takeuchi equivalent $k$-coalgebras;
\item There is a $k$-linear equivalence between the categories of pseudocompact left modules over $\widehat{A} = (A^\circ)^*$ and $\widehat{B} = (B^\circ)^*$.
\end{enumerate}
Then \textnormal{(a)$\implies$(b)$\iff$(c)}.
\end{proposition}

\begin{proof}
The equivalence (b)$\iff$(c) is a direct consequence of the $k$-linear duality~\cite[Theorem~4.3]{simson} between the categories of left $A^\circ$-comodules and pseudocompact left modules over $(A^\circ)^* \cong \widehat{A}$. Thus it suffices to show that (a)$\implies$(c).

Recall~\cite[Proposition~18.33]{lam:lectures} that two rings $R$ and $S$ are Morita equivalent if and only if $S \cong e \M_n(R) e$ for some integer $n \geq 1$ and a full idempotent $e \in \M_n(R)$. Assuming~(a) holds, we may thus write
\begin{equation}\label{eq:corner}
B \cong e \M_n(A) e
\end{equation}
for $n \geq 1$ and a full idempotent $e \in \M_n(A)$, where the above is a $k$-algebra isomorphism by $k$-linearity of the Morita equivalence. For each open ideal $I \in \F(B)$ in the cofinite topology on $B$, there is a corresponding $J \in \F(A)$ such that, under~\eqref{eq:corner}, both $I \cong e \M_n(J) e$ and $B/I \cong e \M_n(A/J) e$ (where we view $\M_n(A/J)$ as a bimodule over $\M_n(A)$ via its natural surjection onto $\M_n(A)/\M_n(J) \cong \M_n(A/J)$). Thus if we pass to the completion, we obtain an isomorphism of topological algebras
\begin{align*}
\widehat{B} &= \invlim\nolimits_{I \in \F(B)} B/I \\
&\cong \invlim\nolimits_{J \in \F(A)} e \M_n(A/J) e \\
&\cong e \M_n(\widehat{A}) e,
\end{align*}
where we now view $\M_n(\widehat{A})$ as a bimodule over $\M_n(A)$ via the algebra homomorphism given by completion $\M_n(A) \to \widehat{\M_n(A)} \cong \M_n(\widehat{A})$.
Letting $e_{ij} \in \M_n(\widehat{A})$ denote the matrix units, we define pseudocompact bimodules
\begin{align*}
{}_{\widehat{A}} U_{\widehat{B}} &= e_{11} \M_n(\widehat{A})e, \\
{}_{\widehat{B}} V_{\widehat{A}} &= e \M_n(\widehat{A}) e_{11}.
\end{align*}
It is then straightforward to verify using the completed tensor product of pseudocompact bimodules (see~\cite[Section~2]{brumer:pseudocompact} or~\cite[Section~4]{vandenbergh:blowup}) we have $U \ctens_{\widehat{B}} V \cong \widehat{A}$ and $V \ctens_{\widehat{A}} U \cong \widehat{B}$, so that the $k$-linear functors
\begin{align*}
U \otimes_{\widehat{B}} -  \colon \widehat{B}\lPC &\to \widehat{A}\lPC,\\
V \otimes_{\widehat{A}} - \colon \widehat{A}\lPC &\to \widehat{B}\lPC
\end{align*}
yield a $k$-linear equivalence of categories. Thus~(c) holds as desired.
\end{proof}

\section{Affine algebras finite over their center}\label{sec:Azumaya}

In this final section, we specialize to the case of affine algebras $R$ that are module-finite over their center $Z(R)$, with the goal of illustrating how the spectrum of $Z(R)$ is reflected in the structure of $R^\circ$. 
First we show that a large part of $R^\circ$ is controlled by an open subscheme of $\Spec Z(R)$ defined in terms of its Azumaya locus.
We then apply this method explicitly to the case of the quantum plane~\cite[1.2]{manin:quantum} at a root of unity. 
Finally, we discuss how a similar analysis can be applied to other algebras if one has a sufficiently good understanding of the local structure of $R$ relative to $\Max Z(R)$.

In order to clarify the geometric picture, we will \textbf{assume $k$ is algebraically closed }throughout this section.
Suppose that $R$ is an affine $k$-algebra that is module-finite over its center $Z(R)$. It follows (by the Artin-Tate Lemma) that $Z(R)$ is affine and consequently noetherian. So $R$ is an affine noetherian PI algebra and therefore is fully RFD by Proposition~\ref{prop:FBN RFD}.

\subsection{The Azumaya locus in the dual coalgebra}\label{sub:azumaya}

To provide a clearer picture of the finite dual, we wish to incorporate information about the representation theory of an algebra.
%
We will do so for an algebra $R$ under the following hypothesis:
\begin{itemize}
\item[(\textbf{H})] $R$ is a prime affine $k$-algebra that is module-finite over its center $Z(R)$.
\end{itemize}
If $R$ satisfies~(\textbf{H}), it is known~\cite[10.2]{mcconnellrobson} that the prime spectrum of $R$ is closely related to $\Spec Z(R)$. 
Our goal here is to similarly relate the quantized spectrum $R^\circ$ to the maximal spectrum of its center. Specifically, we will describe how the Azumaya locus in $\Max Z(R)$ appears within $R^\circ$.

We recall some of the theory of the Azumaya locus from~\cite[Section~3]{browngoodearl:pi} and~\cite[Section~III.1]{browngoodearl:quantum}. Suppose that $R$ satisfies hypothesis~(\textbf{H}).
Let $d$ denote the PI degree~\cite[13.6]{mcconnellrobson} of $R$. Then $d$ is equal to the maximal $k$-dimension of all simple $R$-modules. The \emph{Azumaya locus} is the subset of the maximal spectrum $\Max(Z(R))$ that can be characterized in the following equivalent ways:
\begin{align*}
\Azu(R) &= \{\m \in \Max(Z(R)) \mid R_\m \mbox{ is an Azumaya algebra over } Z(R)_\m\} \\
&= \{\m \in \Max(Z(R)) \mid \m = Z(R) \cap \ann(V) \mbox{ for a simple } {}_R V,\  \dim_k(V) = d\} \\
&= \{\m \in \Max(Z(R)) \mid R\m \mbox{ is a maximal ideal of } R\} \\
&= \{\m \in \Max(Z(R)) \mid R/R\m \cong \M_d(k)\}.
\end{align*}
This is a nonempty open (and therefore dense) subset of $\Max Z(R)$. 
Under an additional mild hypothesis (which holds automatically in the case $Z(R)$ is normal), it is shown in~\cite{brownyakimov:azumaya} that $\Azu(R)$ is the complement in $\Max Z(R)$ of a \emph{discriminant} ideal $D_{n^2}(R/Z(R), \operatorname{tr}) \subseteq Z(R)$.

A prime ideal $P$ of $R$ is \emph{regular} if the PI~degree of $R/P$ is also equal to $d$.
Then a maximal ideal $M$ of $R$ is regular if and only if $M = R\m$ for some $\m \in \Azu(R)$. 
Lemma~\ref{lem:azumaya point} below will provide a link between $R^\circ$ and the points in the Azumaya locus of $R$.

\begin{lemma}\label{lem:matrix azumaya}
Let $(C,\m)$ be a commutative complete local ring whose residue field is algebraically closed. Then every Azumaya algebra over $C$ is isomorphic to a matrix ring over $C$. 
\end{lemma}

\begin{proof}
By completeness of $C$, the morphism of Brauer groups $\Br(C) \to \Br(C/\m)$ induced by the surjection $C \twoheadrightarrow C/\m$ is injective~\cite[Corollary~6.2]{auslandergoldman:brauer}. The Brauer group of the algebraically closed field $C/\m$ is trivial, so $\Br(C)$ is also trivial. It follows~\cite[Proposition~5.3]{auslandergoldman:brauer} that every Azumaya algebra over $C$ is the endomorphism ring of a finitely generated projective module $P_C \neq 0$. Because $C$ is local, $P$ is in fact free, so that the conclusion follows.
\end{proof}

\begin{lemma}\label{lem:azumaya point}
Let $R$ be an algebra over an algebraically closed field $k$ satisfying~\textnormal{(\textbf{H})}, and let $d$ be the PI~degree of $R$. Let $M$ be a regular maximal ideal of $R$ and set $\m = M \cap Z(R) \in \Azu(R)$. Then we have
\[
\dirlim \, (R/M^i)^* \cong \M^d \otimes \Dist(\Spec Z(R),\m).
\]
\end{lemma}

\begin{proof}
Denote $C = Z(R)$.
The algebra $R/M^i \cong (R/M^i)_\m \cong R_\m/(M_\m)^i$ is a homomorphic image of $R_\m$, which is Azumaya over $C_\m$. Note that the completion of $R_\m$ in the cofinite topology coincides with its completion in the $\m$-adic topology (since its unique maximal ideal is $M_\m = R_\m \m$). Because $C_\m$ is noetherian, this completion~\cite[Theorem~7.2]{eisenbud} is given by
\[
\widehat{R_\m} \cong \widehat{C_\m} \otimes_{C_\m} R_\m. 
\]
Because the Azumaya property is stable under extension of scalars~\cite[Corollary~1.6]{auslandergoldman:brauer}, it follows that $\widehat{R_\m}$ is Azumaya over the complete local ring $\widehat{C_\m}$, whose residue field $k$ is algebraically closed. Thus Lemma~\ref{lem:matrix azumaya} implies that $\widehat{R_\m}$ is a matrix algebra over $\widehat{C_\m}$. Since $\widehat{R_\m}/\widehat{R_\m}\m \cong R/R\m \cong \M_d(k)$, we must in fact have 
\[
\widehat{R_\m} \cong \M_d(\widehat{C_\m}) .
\]
Thus the inversely directed system of algebras
\[
R/M^i \cong R_m/M_\m^i \cong \widehat{R_m}/(\widehat{R_\m}\m)^i 
\]
is isomorphic to the system of algebras 
\[
\M_d(\widehat{C_\m}/\widehat{C_\m}\m^i) \cong \M_d(C_\m/\m_\m^i) \cong \M_d(k) \otimes (C_\m/\m_\m^i).
\]
The claim now follows because $\Dist(C, \m) \cong C_\m^\circ \cong \dirlim (C_\m/\m_\m^i)^*$.
\end{proof}

We are now prepared to prove the key result describing the relationship between the Azumaya locus and $R^\circ$. To facilitate its statement, we will define an open subscheme of $\Spec Z(R)$ whose closed points are exactly $\Azu(R)$. 
Denote the intersection of all non-regular maximal ideals of $R$ by
\[
N = \bigcap \{M \in \Max(R) \mid M \mbox{ is not regular}\}.
\]
Then~\cite[Lemma~III.1.2]{browngoodearl:quantum} a prime ideal of $R$ is non-regular if and only if it contains $N$.
Combining this observation with the characterizations of the Azumaya locus as well as the fact~\cite[Lemma~III.1.5]{browngoodearl:quantum} that the contraction map from the prime spectrum of $R$ to $Z(R)$ is closed, it follows that a maximal ideal $\m$ of $Z(R)$ lies outside of the Azumaya locus if and only if it contains the ideal 
\[
I_R = N \cap Z(R).
\]
Let $U_R = D(I_R)$ be the open subscheme of $\Spec Z(R)$ that is the non-vanishing locus of $I_R$. We will call this the \emph{Azumaya subscheme} of $\Spec Z(R)$.
By the discussion above, its set of closed points is exactly equal to $\Azu(R)$.  

\begin{theorem}\label{thm:azumaya}
Let $R$ be an algebra over an algebraically closed field $k$ satisfying~\textnormal{(\textbf{H})}. Let $d$ denote the PI degree of $R$, let $N \unlhd R$ denote the intersection of all non-regular maximal ideals in $R$, and let $U_R$ be the  Azumaya subscheme of $\Spec Z(R)$ described above. 
Then there is an isomorphism of coalgebras
\[
R^\circ \cong \dirlim_{i \geq 1} \, (R/N^i)^\circ  \oplus (\M^d \otimes \Dist(U_R)).
\]
\end{theorem}

\begin{proof}
Denote $C = Z(R)$, a commutative affine algebra over which $R$ is a finitely generated module.
Let $I \in \F(R)$. The image of the composite $C \to R \to R/I$ is a finite-dimensional commutative algebra and thus has the form 
\[
C_I := C/(C \cap I) \cong \prod_{i=1}^s C/\m_i^{e_i}
\] 
for distinct maximal ideals $\m_i$ of $C$. Up to reordering, suppose that $\m_1,\dots,\m_r \notin \Azu(R)$ and $\m_{r+1}, \dots, \m_s \in \Azu(R)$. Then we may fix idempotents $e_I \in C_I$ and $f_I = 1-e_I \in C_I$ so that via the isomorphism above we have
\begin{equation}\label{eq:idempotents}
e_I C_I \cong \prod_{i=1}^r C/\m_i^{e_i} \qquad \mbox{and} \qquad
f_I C_I \cong \prod_{i=r+1}^s C/\m_i^{e_i}.
\end{equation}
Since $e_I$ and $f_I$ are central in $\overline{R} := R/I$, we have $\overline{R} = e_I(R/I) \oplus f_I(R/I)$. 

By construction, the maximal ideals of $f_I \overline{R}$ correspond to regular maximal ideals of $R$, while those of $e_I \overline{R}$ correspond to non-regular maximal ideals of $R$. 
Furthermore, if $I \subseteq I'$ are ideals of finite codimension in $R$, then the image of $e_I \in R/I$ under $R/I \twoheadrightarrow R/I'$ is $e_{I'}$, and similarly $f_I$ maps to $f_{I'}$. Thus the decomposition $(R/I)^* \cong (e_I R/I)^* \oplus (f_I R/I)^*$ is compatible with the directed limit taken over the $I \in \F(R)$, from which it follows that 
\[
R^\circ \cong \dirlim_{I \in \F(R)} (e_I R/I)^* \oplus \dirlim_{I \in \F(R)} (f_I R/I)^*.
\]
Let $\F_n \subseteq \F(R)$ denote those ideals of finite codimension that are contained only in non-regular maximal ideals of $R$, and let $\F_r \subseteq \F(R)$ be those ideals that are contained only in regular maximal ideals. Then the decomposition above amounts to
\begin{equation}\label{eq:regular decomp}
R^\circ \cong \dirlim_{I \in \F_n} (R/I)^* \oplus \dirlim_{I \in \F_r} (R/I)^*.
\end{equation}

First we will prove that $\dirlim_{I \in \F_n} (R/I)^* \cong \dirlim_i (R/N^i)^\circ$. To do so, it is enough to prove that 
\[
\F_n = \{I \in \F(R) \mid N^i \subseteq I \mbox{ for some } i \geq 1\}.
\]
First suppose that $I \in \F_n$, and let $M_1,\dots,M_n \subseteq R$ be the maximal ideals above $I$, so that $N \subseteq M_1 \cap \cdots \cap M_n$. Then $(M_1 \cap \cdots M_n)/I$ is the Jacobson radical of $R/I$ which is nilpotent, say of order $i$. This means that $N^i \subseteq (M_1 \cap \cdots \cap M_n)^i \subseteq I$.
Conversely, suppose that $I \in \F(R)$ with some $N^i \subseteq I$. Then every maximal ideal containing $I$ also contains $N$ and therefore is not regular. So $I \in \F_n$, establishing the desired equality.

Now we describe $\dirlim_{I \in \F_r} (R/I)^*$.
By nilpotence of the Jacobson radical of each $R/I$, this is the same as computing the direct limit of the $(R/I)^*$ where $I$ is a product of regular maximal ideals $M = R\m$ for some $\m \in \Azu(R)$. 
Given such maximal ideals $M_i = R\m_i$, note that their products commute in the monoid of ideals of $R$ (as $\m_i \subseteq C$ are central). Then if $I$ is a product of powers of distinct such $M_i$, we have 
\[
R/I = R/(M_1^{e_1} \cdots M_r^{e_r}) \cong \bigoplus R/M_i^{e_i},
\]
so that
\[
(R/I)^* \cong \bigoplus (R/M_i^{e_i})^*.
\]
Taking the direct limit over all $I$ is the same as taking the limit over all possible products of maximal ideals, which combines with Lemma~\ref{lem:azumaya point} to result in
\begin{align*}
\dirlim_{I \in \F_r} (R/I)^* &\cong \bigoplus_{\m \in \Azu(R)} \dirlim_i \, (R/R\m^i)^* \\
&\cong \bigoplus_{\m \in \Azu(R)} \M^d \otimes \Dist(\Spec C,\m) \\
&\cong \M^d \otimes \bigoplus_{\m \in \Azu(R)} \Dist(\Spec C,\m).
\end{align*}
Finally, taking into account that $\Azu(R)$ consists precisely of the closed points of the open subscheme $U = U_R$ of $\Spec C$ and that the distributions based at the point $\m$ are independent of the open subscheme in which we compute them, we have
\[
\bigoplus_{\m \in \Azu(R)} \Dist(\Spec C,\m) = \bigoplus_{\m \in U_0} \Dist(U,\m) \cong \Dist(U)
\]
thanks to Proposition~\ref{prop:distribution sum}. So $\dirlim_{I \in \F_r} (R/I)^* \cong \M^d \otimes \Dist(U_R)$, and the conclusion now follows from~\eqref{eq:regular decomp}.
\end{proof}

How should we interpret isomorphism of Theorem~\ref{thm:azumaya}? The canonical surjection $R \twoheadrightarrow R/N$ yields an embedding $(R/N)^\circ \hookrightarrow R^\circ$, which we can view as the inclusion of a closed subspace~\cite[Section~4]{smith:subspaces}. The direct limit $\dirlim (R/N^i)^\circ$ thus represents the underlying discrete part of a formal neighborhood of this closed subspace within $R^\circ$. 
The complementary summand is the underlying discrete part of a complementary open subspace. Thus it can be imagined as a ``direct product''~\cite[1.1]{manin:quantum} of a qu$d$it (where $d$ is determined by the representation theory of $R$) with the Azumaya locus (an open subspace of the maximal spectrum of the center).
In light of the decomposition of $\corad R$ from~\eqref{eq:qubit part}, this summand contains the qu$d$its within the spectrum of $R$, while $(R/N)^\circ$ contains the qu-$n$its for all values $n < d$. 

This specializes nicely to the case of an Azumaya algebra. By the Artin-Procesi Theorem~\cite[III.1.4]{browngoodearl:quantum}, this is the case where $R/N = 0$, so that summand corresponding to the formal neighborhood of that closed subspace is trivial. In this case the Azumaya locus is the entire maximal spectrum $\Azu(R) = \Max Z(R)$, so that the Azumaya subscheme is $U_R = \Spec Z(R)$. So for Azumaya algebras the center governs the entire structure of the dual coalgebra in the following way. This recovers a version of~\cite[Example~2]{lebruyn:branes}.

\begin{corollary}
Suppose that $R$ is an algebra over an algebraically closed field $k$ satisfying~\textnormal{(\textbf{H})}, and let $d$ be its PI~degree. If $R$ is Azumaya over its center $Z(R)$, then
\[
R^\circ \cong \M^d \otimes \Dist(\Spec Z(R)).
\]
\end{corollary}


\subsection{Picturing the quantum plane}\label{sub:quantum plane}

In this subsection we will analyze the dual coalgebra of the quantum plane at a root of unity, using the Azumaya locus technique above. From the perspective set out in Sections~\ref{sec:discrete} and~\ref{sec:quantized Max} above, this will give us a glimpse of the quantized set of closed points in the quantum plane. 

We continue to assume that our field $k$ is algebraically closed. In the discussion below, we will alternatively view the affine plane over $k$ as the classical algebraic variety $k^2$ and as the scheme $\A^2_k = \Spec k[x,y]$, depending on the best context for a particular observation. We trust that this fluctuating perspective can be reasonably navigated by the careful reader. 

Let $q \in k^*$. Recall that the (algebra of functions on the) \emph{quantum plane} is the affine domain
\[
\O_q(k^2) = k_q[x,y] = k \langle x, y \mid yx = q  xy \rangle.
\]
Suppose  that $q$ is a primitive $n$th root of unity, so that if $k$ has characteristic $p > 0$ then $n$ is relatively prime to $p$. A number of facts stated without justification below can be found in~\cite[Section~7]{deconciniprocesi:quantum} and~\cite[Examples~III.1.2(3), III.1.4]{browngoodearl:quantum}. 
The center of the quantum plane is given by
\[
Z(\O_q(k^2)) = k[x^n, y^n].
\]
The algebra $\O_q(k^2)$ is module-finite over its center and satisfies~(\textbf{H}). Furthermore, the PI~degree of $\O_q(k^2)$ is equal to~$n$.

While the center $Z(\O_q(k^2)) = k[x^n, y^n]$ is abstractly isomorphic to the coordinate ring of the affine plane $\A^2_k$, we will not view it in this way because we wish to identify the subalgebras $k[x]$ and $k[y]$ of $\O_q(k^2)$ with the coordinate rings of the axes within the ordinary affine plane.
Instead, we will view the center as the coordinate ring of a quotient of the affine plane by a group action.

Let $G_q = \langle q \rangle^2 \subseteq (k^*)^2$ denote the subgroup of the $2$-torus whose coordinates are powers of $q$. This is a finite group with $G_q \cong (\Z/n\Z)^2$. The usual action of the 2-torus on $k^2$ by coordinate scaling restricts to an action of $G_q$ on the plane:
\begin{equation}\label{eq:group action}
(q^i, q^j) \cdot (\alpha, \beta) = (q^i \alpha, q^j \beta).
\end{equation}
This corresponds to an action of $G_q$ on the coordinate ring $k[x,y] = \O(k^2)$ given by $(q^i,q^j) \cdot f(x,y) = f(q^i x, q^j y)$, which has fixed subalgebra 
\[
k[x,y]^{G_q} = k[x^n,y^n] \cong Z(\O_q(k^2)).
\] 
For this reason we view the center as the coordinate ring of the quotient scheme
\[
\Spec k[x^n,y^n] = \A^2_k/G_q.
\]

To explicitly describe the coradical of $\O_q(k^2)^\circ$ as in~\eqref{eq:qubit part} amounts to enumerating the finite-dimensional irreducible representations of the quantum plane. 
These have been parametrized, for instance, in~\cite[Section~7]{deconciniprocesi:quantum}, but we have been unable to find an explicit list of these irreps with the exception of~\cite[Section~8]{reichwalton:sklyanin}, which studies the case where $k = \C$ and $q = -1$. Thus we pause to record the following.

\begin{lemma}\label{lem:irreps}
Let $k$ be an algebraically closed field, and let $q \in k^\times$ be a primitive $n$th root of unity. The irreducible representations $V(\alpha, \beta)$ of $R = \O_q(k^2)$ are parametrized by $(\alpha, \beta) \in k^2$ and are given by a homomorphism $R \to \End_k(V(\alpha,\beta))$ as follows:
\begin{enumerate}
\item For $\alpha\beta = 0$, we obtain pairwise inequivalent $1$-dimensional representations, with $R \to k$ given by $x \mapsto \alpha$ and $y \mapsto \beta$. 
\item For $\alpha\beta \neq 0$, there are irreducible representations of the form $V(\alpha,\beta) = k^n$ with $x \mapsto \alpha D$ and $y \mapsto \beta P$, where
\[
D = \begin{pmatrix}
1 & & & \\
 & q & & \\
 & & \ddots & \\
 & & & q^{n-1}
 \end{pmatrix}
 \quad and \quad
P = \begin{pmatrix}
0 & 1 & \cdots & 0 \\
\vdots & \vdots & \ddots & \vdots \\
0 & 0 & \cdots & 1 \\
1 & 0 & \cdots & 0
\end{pmatrix}.
\]
Two of these representations satisfy $V(\alpha, \beta) \cong V(\alpha', \beta')$ if and only if $(\alpha, \beta)$ and $(\alpha',\beta')$ lie in the same orbit of the $G_q$-action~\eqref{eq:group action} the plane.
\end{enumerate}
\end{lemma}

\begin{proof}
Let $M$ be a maximal ideal of $R$, and denote $\m = M \cap k[x^n, y^n]$. Then $\m = (x^n - c, y^n - d)$ for some $(c,d) \in k^2$. Because $k$ is algebraically closed, we may fix $(\alpha,\beta) \in k^2$ such that $(c,d) = (\alpha^n, \beta^n)$, so that
\[
\m = (x^n - \alpha^n, y^n - \beta^n). 
\]
Let $S = R/R\m$, which is a finite-dimensional algebra with $k$-basis $\{x^iy^j \mid 0 \leq i, j \leq n-1\}$. 

(1) Suppose $\alpha\beta = 0$. Without loss of generality, we may assume that $\beta = 0$ (the case $\alpha = 0$ being symmetric). Then $\m = (x^n - \alpha^n, y^n)$. Note that the image of $y$ lies in the Jacobson radical of $S$.
Thus $S$ has the same irreducible representations as $S/(y) \cong k[x]/(x^n - \alpha^n)$, whose irreps are all $1$-dimensional and given by evaluating $x$ at the roots of $x^n - \alpha^n = x^n - c$, one of which is $\alpha$.
Thus we obtain the representation $S \to k$ given by $x \mapsto \alpha$ and $y \mapsto \beta$ whenever $\alpha\beta = 0$.

(2) Now suppose that $\alpha\beta \neq 0$. In this case $\m = (x^n-\alpha^n, y^n-\beta^n)$ lies in the Azumaya locus $\Azu(R)$, so that it lies below the single maximal ideal $M = R\m$.
One may readily verify that for $D$ and $P$ as in the statement above, the matrices $X = \alpha D$ and $Y = \beta P$ satisfy $YX = qXY$ as well as $X^n = \alpha^n I$ and $Y^n = \beta^n I$. Thus for $V(\alpha,\beta) = k^n$ we have an algebra homomorphism $\phi \colon R \to \End_k(V(\alpha,\beta))$ given by $x \mapsto X$ and $y \mapsto Y$ whose kernel equals $M$. To see that this representation is irreducible, we only need to verify that $\phi$ is surjective. This follows easily from the fact that the regular maximal ideal $M$ has codimension $n^2$. (It can also be verified explicitly by noting that $\phi(\alpha^{-1}x) = D$ and $\phi(\beta^{-1}y) = P$ generate the algebra $\End_k(V)$, an argument that is facilitated by noticing that $\frac{1}{n}(I + D + D^2 + \cdots + D^{n-1}) = E_{11}$ is a standard matrix unit.)

Every algebra automorphism of $S = R/M \cong \End_k(k^n)$ is inner and thus is induced by an intertwiner. So we have $V(\alpha, \beta) \cong V(\alpha', \beta')$ if and only if the kernels of the representations are the same regular maximal ideals, if and only if the central maximal ideals $(x^n - \alpha^n, y^n-\beta^n) = (x^n - (\alpha')^n, y^n-(\beta')^n)$ are equal. This occurs exactly when $\alpha^n = (\alpha')^n$ and $\beta^n = (\beta')^n$, which is equivalent to saying that the points $(\alpha,\beta)$ and $(\alpha', \beta')$ are in the same $G_q$-orbit.
\end{proof}

Thus we can describe the structure of the coradical quite precisely in the following way. Let $W = V(xy) \subseteq k^2$ be the union of the coordinate axes, and let $Y_q \subseteq k^2/G_q$ be the open subset of the quotient space that is the image of the $G_q$-invariant open complement $k^2 \setminus W \subseteq k^2$. Then Lemma~\ref{lem:irreps} yields:
\begin{align}\label{eq:plane coradical}
\begin{split}
\corad \O_q(k^2) &\cong \bigoplus_{W} k \oplus \bigoplus_{Y_q} \M^n \\
&\cong kW \oplus (\M^n \otimes kY_q).
\end{split}
\end{align}

This hints at a decomposition of the full finite dual, which can be rigorously demonstrated using the method of Subsection~\ref{sub:azumaya}. Let $N$ be the intersection of the non-regular maximal ideals in $\O_q(k^2)$. From Lemma~\ref{lem:irreps}(1) we have
\[
N = (x) \cap (y) = (xy),
\]
whose intersection with the center is $(xy) \cap k[x^n, y^n] = (x^n y^n)$.
Let $V = V(xy)$ be the closed subscheme of $\A^2_k$ that is the union of the coordinate axes. 
If we view $\Spec Z(R) = \A^2_k/G_q$ as discussed above, then the Azumaya subscheme $U_q = D(xy)/G_q$ is the open subscheme that is the image of the $G_q$-invariant complement of $V$. Theorem~\ref{thm:azumaya} now gives us
\begin{equation}\label{eq:quantum plane}
\O_q(k^2)^\circ \cong \dirlim_{i \geq 1} \, (\O_q(k^2)/(xy)^i)^\circ \oplus (\M^n \otimes \Dist U_q).
\end{equation}
Since the varieties $W$ and $Y_q$ precisely correspond to the closed points of the schemes $V$ and $U_q$, this decomposition restricts to~\eqref{eq:plane coradical} when restricting to the coradical.

Note that the algebra $\O_q(k^2)/N \cong k[x,y]/(xy)$ is independent of $q$, commutative, and isomorphic to the coordinate ring of $V = V(xy)$. In particular,
\[
(\O_q(k^2)/N)^\circ \cong \Dist V.
\]
Thus the first summand in~\eqref{eq:quantum plane} is a formal neighborhood of the classical union of the coordinate axes. But even though the summand of the coradical~\eqref{eq:plane coradical} corresponding to the axes is classical (i.e., cocommutative), its formal neighorhood has a truly quantum (non-cocommutative) nature. For instance, its first-order neighborhood $(\O_q(k^2)/N^2)^\circ$ contains the dual of the 4-dimensional algebra $\O_q(k^2)/(x^2, y^2)$ which is not commutative if $q \neq 1$.

\begin{figure}
\includegraphics{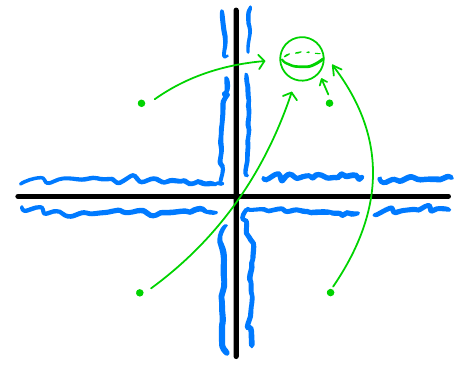}
\caption{Depiction of $\O_q(\C^2)^\circ$ at $q = -1$}
\label{fig:quantum plane}
\end{figure}

Figure~\ref{fig:quantum plane} illustrates the real part of the complex quantum plane $\O_q(\C^2)^\circ$ in the case where $q = -1$, so that $n = 2$. The axes are classical sets of points, but with quantum formal neighborhoods. Meanwhile, the points off of the axes have orbits under the action of the group $G_{-1} = \{\pm 1\}^2$ of order~$4$. The points in each orbit become identified in the quotient by the group action and then ``replaced'' by the qubit $\M^2$.

\subsection{Further discussion}
\label{sub:orders}

In this final section we discuss how a similar analysis can in principle be carried out for other affine algebras that are module-finite over their center. Outside of the Azumaya locus, the finite dual can be described locally relative to the center with the following method. (This generalizes Proposition~\ref{prop:distribution sum}.)

%
%
%
%

\begin{proposition}\label{prop:local structure}
Let $R$ be an affine $k$-algebra that is module-finite over its center $Z = Z(R)$. Then 
\[
R^\circ \cong \bigoplus_{\m \in \Max Z} \dirlim \, (R/\m^i R)^* \cong \bigoplus_{\m \in \Max Z} (R_\m)^\circ.
\]
\end{proposition}

\begin{proof}
This proof is similar in spirit to that of Theorem~\ref{thm:azumaya}. Consider the ideals of $R$ of the form $JR$ where $J \in \F(Z)$. These have finite codimension in $R$ since $R/JR \cong Z/J \otimes_Z R$ is module-finite over $Z/J$, and they form a cofinal subset in $\F(R)$ since each $I \in \F(R)$ has $J := Z \cap I \in \F(Z)$ and $JR \subseteq I$. 

As discussed in the proof of Theorem~\ref{thm:azumaya}, each $J \in \F(Z)$ has the form $J = \prod_{j=1}^s \m_j^{e_j}$ for some $\m_j \in \Max Z$ and integers $e_j \geq 0$. Thus our cofinal family of ideals in $\F(R)$ have the form
\[
JR = (\m_1^{e_1} \cdots \m_s^{e_s}) R = (\m_1^{e_1} R) \cdots (\m_s^{e_s} R).
\]
From the Chinese Remainder Theorem we obtain an isomorphism of coalgebras
\[
(R/ \m_1^{e_1} \cdots \m_s^{e_s} R)^* \cong \bigoplus_{j=1}^s (R/\m_j^{e_j} R)^*.
\]
Passing to the direct limit over all $J \in \F(Z)$, which amounts to all possible (finite) products of powers of maximal ideals in $Z$, we obtain one of our desired isomorphisms:
\[
R^\circ \cong \dirlim_{J \in \F(Z)} (R/JR)^* \cong \bigoplus_{\m \in \Max Z} \dirlim_i \, (R/\m^i R)^*.
\]
The second isomorphism of the statement above follows from $(R_\m)^\circ \cong \dirlim (R/\m^i R)^*$, which is readily verified for any $\m \in \Max Z$.
\end{proof}

By the duality~\cite[Theorem~3.6]{simson} between coalgebras and pseudocompact algebras, to understand $(R_\m)^\circ$ is equivalent to understanding the $\m$-adic completion $\widehat{R_\m} \cong (R_\m)^{\circ *}$. Thus $R^\circ$ can in principle be described as above if one has a good understanding of the completed local structure of $R$ relative to the maximal spectrum of its center. 

The complete local structure of $R$ is in turn determined by its \e tale local structure. The \e tale local structure of orders $R$ over schemes has been of interest for many years. 
An important case investigated by Le Bruyn~\cite{lebruyn:local} is the local structure of Cayley-Hamilton algebras. In the case where $R$ is smooth in the sense of~\cite[Section~4]{procesi:cayley}, he showed that $R$ has only finitely many \e tale-local isomorphism classes. 
Much more about these orders can be found in~\cite{lebruyn:ncg}.

Another special case in which the \e tale local structure has received much attention is that of orders over surfaces~\cite{artin:surfaces, chaningalls:minimal, chanhackingingalls:singularities}. We will close with some remarks on how the known structure theory in this case relates to the decomposition of Proposition~\ref{prop:local structure}.

Let $X$ be a smooth affine surface over~$k$, and let $R$ be a maximal order over $Z = k[X]$. Then $R$ is ramified on the  curve $D$ which is the zero locus in $X$ of the reduced discriminant $I = \sqrt{d(R)}$, where $d(R)$ is the classical discriminant~\cite[\S 10]{reiner}. 
As described in~\cite[Proposition~1.6]{ingalls:quantizable}, for any closed point $x$ of $X$, we are in one of the following three cases:
\begin{enumerate}[label=(\roman*)]
\item $x \notin D$, and $R$ is \e tale locally isomorphic to a matrix algebra at $x$; 
\item $x$ is a smooth point of $D$, and $R$ is \e tale locally isomorphic at~$x$ to a hereditary order (whose structure is described in~\cite[(39.14)]{reiner}); 
\item $x$ is a singular point of $D$.
\end{enumerate}
All three of these cases appear in the example of the quantum plane $\O_q(k^2)$ above. This is an order over the surface $X = \A^2_k/G_q$, with notation as in the previous subsection. Here $D$ is the image in $X$ of the $G_q$-invariant closed subscheme $W = V(xy)$ of $\A^2_k$, which contains the image of the $G_q$-invariant origin.
Then case~(i) occurs on the complement of $D$ (which is the Azumaya locus), case~(ii) occurs on $D$ away from the origin, and case~(iii) occurs at the origin. 

For a general smooth surface $X$ with maximal order $R$ over $k[X]$, suppose that a point $x$ falls under case~(iii). If $D$ is a divisor having only normal crossings, it is shown in~\cite[Theorem~1.2]{artin:surfaces} that there are only finitely many possible isomorphism classes for the \e tale local structure at~$x$. On the other hand, under the additional assumption that $R$ has global dimension~2, the possible complete local structures of $R$ have been classified up to Morita equivalence by Artin~\cite{artin:maximal, artin:twodim} and independently by Reiten and Van den Bergh~\cite{reitenvandenbergh:maximal}. This classification has been treated uniformly in~\cite{chaningalls:lowdim} in terms of crossed product orders $k[[u,v]] *_\eta G$ where $G$ is a finite subgroup of $\operatorname{GL}_2(k)$ and $\eta \in H^2(G, k^*)$. Such crossed products are a special case of \emph{twisted tensor product} algebras~\cite{csv:twisted}, whose finite dual coalgebras we will describe elsewhere~\cite{reyes:twisted}.

%

\section{Acknowledgments}

I thank Alex Chirvasitu, Colin Ingalls, Susan Montgomery, Lance Small, Chelsea Walton, and Milen Yakimov for helpful discussions and references to the literature. I also thank So Nakamura and Ari Rosenfield for comments and corrections on an early draft of this paper. 
Finally, I am grateful to the anonymous referee for several thoughtful suggestions that led to a more readable and streamlined version of this paper and spurred the inclusion of Subsection~\ref{sub:orders}.

\bibliography{qmax-v4}

\begin{bibdiv}
\begin{biblist}

\bib{auslandergoldman:brauer}{article}{
      author={Auslander, Maurice},
      author={Goldman, Oscar},
       title={The {B}rauer group of a commutative ring},
        date={1960},
     journal={Trans. Amer. Math. Soc.},
      volume={97},
       pages={367\ndash 409},
         url={https://doi.org/10.2307/1993378},
}

\bib{ananin:representability}{article}{
      author={Anan'in, A.~Z.},
       title={Representability of {N}oetherian finitely generated algebras},
        date={1992},
     journal={Arch. Math. (Basel)},
      volume={59},
      number={1},
       pages={1\ndash 5},
}

\bib{artin:surfaces}{book}{
      author={Artin, M.},
       title={Local structure of maximal orders on surfaces},
      series={Lecture Notes in Math.},
   publisher={Springer, Berlin-New York},
        date={1982},
      volume={917},
}

\bib{artin:maximal}{article}{
      author={Artin, M.},
       title={Maximal orders of global dimension and {K}rull dimension two},
        date={1986},
     journal={Invent. Math.},
      volume={84},
      number={1},
       pages={195\ndash 222},
         url={https://doi.org/10.1007/BF01388739},
}

\bib{artin:twodim}{article}{
      author={Artin, Michael},
       title={Two-dimensional orders of finite representation type},
        date={1987},
        ISSN={0025-2611},
     journal={Manuscripta Math.},
      volume={58},
      number={4},
       pages={445\ndash 471},
         url={https://doi.org/10.1007/BF01277604},
}

\bib{amitsursmall}{article}{
      author={Amitsur, S.~A.},
      author={Small, Lance~W.},
       title={Algebras over infinite fields, revisited},
        date={1996},
     journal={Israel J. Math.},
      volume={96},
       pages={23\ndash 25},
         url={https://doi.org/10.1007/BF02785531},
}

\bib{batchelor:maps}{book}{
      author={Batchelor, Marjorie},
       title={In search of the graded manifold of maps between graded
  manifolds},
      series={Lecture Notes in Phys.},
   publisher={Springer, Berlin},
        date={1988},
      volume={311},
         url={https://doi.org/10.1007/BFb0038539},
}

\bib{batchelor:ncg}{book}{
      author={Batchelor, Marjorie},
       title={Measuring coalgebras, quantum group-like objects and
  noncommutative geometry},
      series={Lecture Notes in Phys.},
   publisher={Springer, Berlin},
        date={1991},
      volume={375},
         url={https://doi.org/10.1007/3-540-53763-5_45},
}

\bib{browngoodearl:quantum}{book}{
      author={Brown, Ken~A.},
      author={Goodearl, Ken~R.},
       title={Lectures on algebraic quantum groups},
      series={Advanced Courses in Mathematics. CRM Barcelona},
   publisher={Birkh\"{a}user Verlag, Basel},
        date={2002},
        ISBN={3-7643-6714-8},
         url={https://doi.org/10.1007/978-3-0348-8205-7},
}

\bib{browngoodearl:pi}{article}{
      author={Brown, K.~A.},
      author={Goodearl, K.~R.},
       title={Homological aspects of {N}oetherian {PI} {H}opf algebras and
  irreducible modules of maximal dimension},
        date={1997},
     journal={J. Algebra},
      volume={198},
      number={1},
       pages={240\ndash 265},
         url={https://doi.org/10.1006/jabr.1997.7109},
}

\bib{bergheunen:obstructions}{article}{
      author={Berg, Benno van~den},
      author={Heunen, Chris},
       title={Extending obstructions to noncommutative functorial spectra},
        date={2014},
     journal={Theory Appl. Categ.},
      volume={29},
      number={17},
       pages={457\ndash 474},
}

\bib{brumer:pseudocompact}{article}{
      author={Brumer, Armand},
       title={Pseudocompact algebras, profinite groups and class formations},
        date={1966},
     journal={J. Algebra},
      volume={4},
       pages={442\ndash 470},
         url={https://doi.org/10.1016/0021-8693(66)90034-2},
}

\bib{brownyakimov:azumaya}{article}{
      author={Brown, K.~A.},
      author={Yakimov, M.~T.},
       title={Azumaya loci and discriminant ideals of {PI} algebras},
        date={2018},
     journal={Adv. Math.},
      volume={340},
       pages={1219\ndash 1255},
         url={https://doi.org/10.1016/j.aim.2018.10.024},
}

\bib{bmr:kochenspecker}{article}{
      author={Ben-Zvi, Michael},
      author={Ma, Alexander},
      author={Reyes, Manuel},
       title={A {K}ochen-{S}pecker theorem for integer matrices and
  noncommutative spectrum functors},
        date={2017},
     journal={J. Algebra},
      volume={491},
       pages={280\ndash 313},
        note={With an appendix by Alexandru Chirvasitu},
}

\bib{chanhackingingalls:singularities}{article}{
      author={Chan, Daniel},
      author={Hacking, Paul},
      author={Ingalls, Colin},
       title={Canonical singularities of orders over surfaces},
        date={2009},
     journal={Proc. Lond. Math. Soc. (3)},
      volume={98},
      number={1},
       pages={83\ndash 115},
         url={https://doi.org/10.1112/plms/pdn025},
}

\bib{chaningalls:minimal}{article}{
      author={Chan, Daniel},
      author={Ingalls, Colin},
       title={The minimal model program for orders over surfaces},
        date={2005},
     journal={Invent. Math.},
      volume={161},
      number={2},
       pages={427\ndash 452},
         url={https://doi.org/10.1007/s00222-005-0438-z},
}

\bib{chaningalls:lowdim}{article}{
      author={Chan, Daniel},
      author={Ingalls, Colin},
       title={Low dimensional orders of finite representation type},
        date={2021},
     journal={Math. Z.},
      volume={297},
      number={3-4},
       pages={1161\ndash 1190},
         url={https://doi.org/10.1007/s00209-020-02552-2},
}

\bib{chinmontgomery:basic}{book}{
      author={Chin, William},
      author={Montgomery, Susan},
       title={Basic coalgebras},
      series={AMS/IP Stud. Adv. Math.},
   publisher={Amer. Math. Soc., Providence, RI},
        date={1997},
      volume={4},
}

\bib{courtneyshulman}{article}{
      author={Courtney, Kristin},
      author={Shulman, Tatiana},
       title={Free products with amalgamation over central
  {$\rm{C}^*$}-subalgebras},
        date={2020},
     journal={Proc. Amer. Math. Soc.},
      volume={148},
      number={2},
       pages={765\ndash 776},
         url={https://doi.org/10.1090/proc/14746},
}

\bib{csv:twisted}{article}{
      author={{\u{C}}ap, Andreas},
      author={Schichl, Hermann},
      author={Van\v{z}ura, Ji\v{r}\'{\i}},
       title={On twisted tensor products of algebras},
        date={1995},
     journal={Comm. Algebra},
      volume={23},
      number={12},
       pages={4701\ndash 4735},
         url={https://doi.org/10.1080/00927879508825496},
}

\bib{deconciniprocesi:quantum}{book}{
      author={De~Concini, C.},
      author={Procesi, C.},
       title={Quantum groups},
      series={Lecture Notes in Math.},
   publisher={Springer, Berlin},
        date={1993},
      volume={1565},
         url={https://doi.org/10.1007/BFb0073466},
}

\bib{demazuregabriel}{book}{
      author={Demazure, Michel},
      author={Gabriel, Peter},
       title={Introduction to {A}lgebraic {G}eometry and {A}lgebraic {G}roups},
      series={North-Holland Mathematics Studies},
   publisher={North-Holland Publishing Co., Amsterdam-New York},
        date={1980},
      volume={39},
        note={Translated from the French by J. Bell},
}

\bib{dnr:hopf}{book}{
      author={D\u{a}sc\u{a}lescu, Sorin},
      author={N\u{a}st\u{a}sescu, Constantin},
      author={Raianu, \c{S}erban},
       title={Hopf {A}lgebras. {A}n introduction},
      series={Monographs and Textbooks in Pure and Applied Mathematics},
   publisher={Marcel Dekker, Inc., New York},
        date={2001},
      volume={235},
}

\bib{eisenbud}{book}{
      author={Eisenbud, David},
       title={Commutative {A}lgebra: {W}ith a view toward algebraic geometry},
      series={Graduate Texts in Mathematics},
   publisher={Springer-Verlag, New York},
        date={1995},
      volume={150},
         url={https://doi.org/10.1007/978-1-4612-5350-1},
}

\bib{goodearlwarfield}{book}{
      author={Goodearl, K.~R.},
      author={Warfield, R.~B., Jr.},
       title={An introduction to noncommutative {N}oetherian rings},
     edition={Second},
      series={London Mathematical Society Student Texts},
   publisher={Cambridge University Press, Cambridge},
        date={2004},
      volume={61},
}

\bib{hormander}{book}{
      author={H\"{o}rmander, Lars},
       title={The {A}nalysis of {L}inear {P}artial {D}ifferential
  {O}perators~{I}},
      series={Grundlehren der mathematischen Wissenschaften [Fundamental
  Principles of Mathematical Sciences]},
   publisher={Springer-Verlag, Berlin},
        date={1983},
      volume={256},
         url={https://doi.org/10.1007/978-3-642-96750-4},
        note={Distribution theory and Fourier analysis},
}

\bib{heunen:classicalfaces}{article}{
      author={Heunen, Chris},
       title={The many classical faces of quantum structures},
        date={2017},
     journal={Entropy},
      volume={19},
      number={4},
       pages={144},
}

\bib{hatipoglulomp:hopf}{article}{
      author={Hatipo\u{g}lu, Can},
      author={Lomp, Christian},
       title={Locally finite representations over {N}oetherian {H}opf
  algebras},
        date={2022},
     journal={Proc. Amer. Math. Soc.},
      volume={150},
      number={5},
       pages={1903\ndash 1923},
         url={https://doi.org/10.1090/proc/15747},
}

\bib{hlfv:enrichment}{article}{
      author={Hyland, Martin},
      author={L\'{o}pez~Franco, Ignacio},
      author={Vasilakopoulou, Christina},
       title={Hopf measuring comonoids and enrichment},
        date={2017},
     journal={Proc. Lond. Math. Soc. (3)},
      volume={115},
      number={5},
       pages={1118\ndash 1148},
         url={https://doi.org/10.1112/plms.12064},
}

\bib{heunenreyes:active}{article}{
      author={Heunen, Chris},
      author={Reyes, Manuel~L.},
       title={Active lattices determine {AW}*-algebras},
        date={2014},
     journal={J. Math. Anal. Appl.},
      volume={416},
      number={1},
       pages={289\ndash 313},
}

\bib{heunenreyes:discretization}{article}{
      author={Heunen, Chris},
      author={Reyes, Manuel~L.},
       title={Discretization of {C}*-algebras},
        date={2017},
     journal={J. Operator Theory},
      volume={77},
      number={1},
       pages={19\ndash 37},
}

\bib{heynemansweedler}{article}{
      author={Heyneman, Robert~G.},
      author={Sweedler, Moss~Eisenberg},
       title={Affine {H}opf algebras. {I}},
        date={1969},
     journal={J. Algebra},
      volume={13},
       pages={192\ndash 241},
         url={https://doi.org/10.1016/0021-8693(69)90071-4},
}

\bib{imr:diagonal}{article}{
      author={Iovanov, Miodrag~C.},
      author={Mesyan, Zachary},
      author={Reyes, Manuel~L.},
       title={Infinite-dimensional diagonalization and semisimplicity},
        date={2016},
     journal={Israel J. Math.},
      volume={215},
      number={2},
       pages={801\ndash 855},
}

\bib{ingalls:quantizable}{article}{
      author={Ingalls, Colin},
       title={Quantizable orders over surfaces},
        date={1998},
     journal={J. Algebra},
      volume={207},
      number={2},
       pages={616\ndash 656},
         url={https://doi.org/10.1006/jabr.1998.7462},
}

\bib{jantzen:groups}{book}{
      author={Jantzen, Jens~Carsten},
       title={Representations of {A}lgebraic {G}roups},
     edition={Second},
      series={Mathematical Surveys and Monographs},
   publisher={American Mathematical Society, Providence, RI},
        date={2003},
      volume={107},
}

\bib{jategaonkar:jacobsons}{article}{
      author={Jategaonkar, Arun~Vinayak},
       title={Jacobson's conjecture and modules over fully bounded {N}oetherian
  rings},
        date={1974},
     journal={J. Algebra},
      volume={30},
       pages={103\ndash 121},
         url={https://doi.org/10.1016/0021-8693(74)90195-1},
}

\bib{johnstone:pointless}{article}{
      author={Johnstone, Peter~T.},
       title={The point of pointless topology},
        date={1983},
     journal={Bull. Amer. Math. Soc. (N.S.)},
      volume={8},
      number={1},
       pages={41\ndash 53},
         url={https://doi.org/10.1090/S0273-0979-1983-15080-2},
}

\bib{kaplansky:operator}{article}{
      author={Kaplansky, Irving},
       title={The structure of certain operator algebras},
        date={1951},
     journal={Trans. Amer. Math. Soc.},
      volume={70},
       pages={219\ndash 255},
}

\bib{kassel:quantum}{book}{
      author={Kassel, Christian},
       title={Quantum {G}roups},
      series={Graduate Texts in Mathematics},
   publisher={Springer-Verlag, New York},
        date={1995},
      volume={155},
}

\bib{kornell:quantumsets}{article}{
      author={Kornell, Andre},
       title={Quantum sets},
        date={2020},
     journal={J. Math. Phys.},
      volume={61},
      number={10},
       pages={102202, 33},
         url={https://doi.org/10.1063/1.5054128},
}

\bib{kornell:discrete}{unpublished}{
      author={Kornell, Andre},
       title={Discrete quantum structures},
        date={2021},
  note={\href{https://arxiv.org/abs/2004.04377v5}{\texttt{arXiv:2004.04377
  [math.OA]}}},
}

\bib{kontsevichsoibelman}{book}{
      author={Kontsevich, M.},
      author={Soibelman, Y.},
       title={Notes on {$A_\infty$}-algebras, {$A_\infty$}-categories and
  non-commutative geometry},
      series={Lecture Notes in Phys.},
   publisher={Springer, Berlin},
        date={2009},
      volume={757},
}

\bib{lam:lectures}{book}{
      author={Lam, T.~Y.},
       title={Lectures on {M}odules and {R}ings},
      series={Graduate Texts in Mathematics},
   publisher={Springer-Verlag, New York},
        date={1999},
      volume={189},
         url={https://doi.org/10.1007/978-1-4612-0525-8},
}

\bib{landsman:aqm}{book}{
      author={Landsman, Nicolaas~P.},
      editor={Greenberger, Daniel},
      editor={Hentschel, Klaus},
      editor={Weinert, Friedel},
       title={Algebraic quantum mechanics},
   publisher={Springer Berlin Heidelberg},
     address={Berlin, Heidelberg},
        date={2009},
        ISBN={978-3-540-70626-7},
         url={https://doi.org/10.1007/978-3-540-70626-7_3},
}

\bib{lebruyn:local}{article}{
      author={Le~Bruyn, Lieven},
       title={Local structure of {S}chelter-{P}rocesi smooth orders},
        date={2000},
        ISSN={0002-9947},
     journal={Trans. Amer. Math. Soc.},
      volume={352},
      number={10},
       pages={4815\ndash 4841},
         url={https://doi.org/10.1090/S0002-9947-00-02567-8},
}

\bib{lebruyn:ncg}{book}{
      author={Le~Bruyn, Lieven},
       title={Noncommutative geometry and {C}ayley-smooth orders},
      series={Pure and Applied Mathematics (Boca Raton)},
   publisher={Chapman \& Hall/CRC, Boca Raton, FL},
        date={2008},
      volume={290},
}

\bib{lebruyn:dual}{unpublished}{
      author={Le~Bruyn, Lieven},
       title={Noncommutative geometry and dual coalgebras},
        date={2008},
        note={\href{https://arxiv.org/abs/0805.2377v1}{\texttt{arXiv:0805.2377
  [math.RA]}}},
}

\bib{lebruyn:branes}{article}{
      author={Le~Bruyn, Lieven},
       title={Representation stacks, {D}-branes and noncommutative geometry},
        date={2012},
     journal={Comm. Algebra},
      volume={40},
      number={10},
       pages={3636\ndash 3651},
}

\bib{manin:quantum}{book}{
      author={Manin, Yu.~I.},
       title={Quantum {G}roups and {N}oncommutative {G}eometry},
   publisher={Universit\'{e} de Montr\'{e}al, Centre de Recherches
  Math\'{e}matiques, Montreal, QC},
        date={1988},
}

\bib{mcconnellrobson}{book}{
      author={McConnell, J.~C.},
      author={Robson, J.~C.},
       title={Noncommutative {N}oetherian {R}ings},
     edition={Revised},
      series={Graduate Studies in Mathematics},
   publisher={American Mathematical Society, Providence, RI},
        date={2001},
      volume={30},
        note={With the cooperation of L. W. Small},
}

\bib{mrv:quantumfunctions}{article}{
      author={Musto, Benjamin},
      author={Reutter, David},
      author={Verdon, Dominic},
       title={A compositional approach to quantum functions},
        date={2018},
     journal={J. Math. Phys.},
      volume={59},
      number={8},
       pages={081706, 42},
         url={https://doi.org/10.1063/1.5020566},
}

\bib{nielsenchuang:quantum}{book}{
      author={Nielsen, Michael~A.},
      author={Chuang, Isaac~L.},
       title={Quantum computation and quantum information},
   publisher={Cambridge University Press, Cambridge},
        date={2000},
        ISBN={0-521-63235-8; 0-521-63503-9},
}

\bib{procesi:cayley}{article}{
      author={Procesi, Claudio},
       title={A formal inverse to the {C}ayley-{H}amilton theorem},
        date={1987},
     journal={J. Algebra},
      volume={107},
      number={1},
       pages={63\ndash 74},
         url={https://doi.org/10.1016/0021-8693(87)90073-1},
}

\bib{radford:hopf}{book}{
      author={Radford, David~E.},
       title={Hopf {A}lgebras},
      series={Series on Knots and Everything},
   publisher={World Scientific Publishing Co. Pte. Ltd., Hackensack, NJ},
        date={2012},
      volume={49},
}

\bib{reiner}{book}{
      author={Reiner, I.},
       title={Maximal {O}rders},
      series={London Mathematical Society Monographs, No. 5},
   publisher={Academic Press [Harcourt Brace Jovanovich, Publishers],
  London-New York},
        date={1975},
}

\bib{reyes:obstructing}{article}{
      author={Reyes, Manuel~L.},
       title={Obstructing extensions of the functor {S}pec to noncommutative
  rings},
        date={2012},
     journal={Israel J. Math.},
      volume={192},
      number={2},
       pages={667\ndash 698},
}

\bib{reyes:sheaves}{book}{
      author={Reyes, Manuel~L.},
       title={Sheaves that fail to represent matrix rings},
      series={Contemp. Math.},
   publisher={Amer. Math. Soc., Providence, RI},
        date={2014},
      volume={609},
         url={https://doi.org/10.1090/conm/609/12129},
}

\bib{reyes:twisted}{unpublished}{
      author={Reyes, Manuel~L.},
       title={Dual coalgebras of twisted tensor products},
        date={2023},
        note={\href{https://arxiv.org/abs/2306.17261}{\texttt{arXiv:2306.17261
  [math.RA]}}},
}

\bib{rump:quantum}{article}{
      author={Rump, Wolfgang},
       title={Symmetric quantum sets and {$L$}-algebras},
        date={2022},
     journal={Int. Math. Res. Not. IMRN},
      number={3},
       pages={1770\ndash 1810},
         url={https://doi.org/10.1093/imrn/rnaa135},
}

\bib{reitenvandenbergh:maximal}{article}{
      author={Reiten, Idun},
      author={Van~den Bergh, Michel},
       title={Two-dimensional tame and maximal orders of finite representation
  type},
        date={1989},
     journal={Mem. Amer. Math. Soc.},
      volume={80},
      number={408},
         url={https://doi.org/10.1090/memo/0408},
}

\bib{reichwalton:sklyanin}{article}{
      author={Reich, Daniel~J.},
      author={Walton, Chelsea},
       title={Explicit representations of 3-dimensional {S}klyanin algebras
  associated to a point of order 2},
        date={2018},
     journal={Involve},
      volume={11},
      number={4},
       pages={585\ndash 608},
         url={https://doi.org/10.2140/involve.2018.11.585},
}

\bib{simson}{article}{
      author={Simson, Daniel},
       title={Coalgebras, comodules, pseudocompact algebras and tame comodule
  type},
        date={2001},
     journal={Colloq. Math.},
      volume={90},
      number={1},
       pages={101\ndash 150},
}

\bib{smith:subspaces}{article}{
      author={Smith, S.~Paul},
       title={Subspaces of non-commutative spaces},
        date={2002},
        ISSN={0002-9947},
     journal={Trans. Amer. Math. Soc.},
      volume={354},
      number={6},
       pages={2131\ndash 2171},
         url={https://doi.org/10.1090/S0002-9947-02-02963-X},
}

\bib{smoktunowicz}{book}{
      author={Smoktunowicz, Agata},
       title={Some results in noncommutative ring theory},
   publisher={Eur. Math. Soc., Z\"{u}rich},
        date={2006},
}

\bib{sweedler:hopf}{book}{
      author={Sweedler, Moss~E.},
       title={Hopf {A}lgebras},
      series={Mathematics Lecture Note Series},
   publisher={W. A. Benjamin, Inc., New York},
        date={1969},
}

\bib{takeuchi:tangent}{article}{
      author={Takeuchi, Mitsuhiro},
       title={Tangent coalgebras and hyperalgebras. {I}},
        date={1974},
     journal={Japan. J. Math.},
      volume={42},
       pages={1\ndash 143},
}

\bib{takeuchi:formal}{article}{
      author={Takeuchi, Mitsuhiro},
       title={Formal schemes over fields},
        date={1977},
     journal={Comm. Algebra},
      volume={5},
      number={14},
       pages={1483\ndash 1528},
}

\bib{takeuchi:morita}{article}{
      author={Takeuchi, Mitsuhiro},
       title={Morita theorems for categories of comodules},
        date={1977},
     journal={J. Fac. Sci. Univ. Tokyo Sect. IA Math.},
      volume={24},
      number={3},
       pages={629\ndash 644},
}

\bib{vasilakopoulu:enrichment}{unpublished}{
      author={Vasilakopoulou, Christina},
       title={Enrichment of categories of algebras and modules},
        date={2012},
        note={\href{https://arxiv.org/abs/1205.6450}{\texttt{arXiv:1205.6450
  [math.CT]}}},
}

\bib{vandenbergh:blowup}{article}{
      author={Van~den Bergh, Michel},
       title={Blowing up of non-commutative smooth surfaces},
        date={2001},
     journal={Mem. Amer. Math. Soc.},
      volume={154},
      number={734},
}

\bib{vanoystaeyen:sheaves}{book}{
      author={Van~Oystaeyen, Freddy M.~J.},
       title={Noncommutative algebraic geometry: a survey of the approach via
  sheaves on noncommutative spaces},
   publisher={Symp. Ring Theory Represent. Theory Organ. Comm., Yamanashi},
        date={2010},
}

\end{biblist}
\end{bibdiv}

\end{document}